\documentclass[11pt]{amsart}
\usepackage{amssymb,amsbsy,amsmath,amsfonts,amssymb,amscd,graphicx}
\usepackage{color}

\newtheorem{theorem}{Theorem}[section]
\newtheorem{corollary}[theorem]{Corollary}
\newtheorem{lemma}[theorem]{Lemma}
\newtheorem{proposition}[theorem]{Proposition}

\numberwithin{equation}{section}

\newcommand{\dx}{\;{\rm d}x}

\newcommand{\ds}{\,{\rm d}s}
\newcommand{\dt}{\,{\rm d}t}

\newcommand{\rd}{{\rm d}}

\newcommand{\dd}{\;{\rm d}}

\def\LL{\mathrm{L}}

\newcommand{\RR}{\mathbb{R}}
\newcommand{\N}{\mathbb{N}}

\def\ee{\mathrm{e}}

\newcommand{\fer}[3]{\mathcal E\left[#1|#2\right]\!#3}
\newcommand{\fe}[2]{\mathcal E[#1]#2}
\newcommand{\I}[3]{\mathcal I[#1|#2]#3}
\newcommand{\EL}[1]{\mbox{\rm \textsf F}[#1]}
\newcommand{\IL}[1]{\mbox{\rm \textsf I}[#1]}
\newcommand{\cf}{{\sl cf.\/} }
\newcommand{\ie}{{\sl i.e.\/}}
\newcommand{\fw}[1]{\mathcal F[{#1}]}
\newcommand{\iw}[1]{\mathcal J[{#1}]}
\newcommand{\dmu}{d\mu}
\newcommand{\dnu}{d\nu}
\newcommand{\m}{\mu}
\newcommand{\iintpolar}{\iint_{(0,\infty)\times\mathbb S^{d-1}}\kern-25pt}
\newcommand{\Cmd}{\mathcal C_{m,d}}
\newcommand{\gopt}{\lambda_{m,d}}
\newcommand{\q}{p}
\newcommand{\Email}[1]{\rm{\sl E-mail\/}: \textsf{#1}}
\newcommand{\overlinetilde}[1]{\widetilde{#1}}

\begin{document}

\title[Asymptotics of the fast diffusion equation]
{Asymptotics of the fast diffusion equation\\
via entropy estimates}

\author[A. Blanchet]{Adrien Blanchet}
\author[M. Bonforte]{Matteo Bonforte}
\author[J. Dolbeault]{Jean Dolbeault}

\author[G. Grillo]{Gabriele Grillo}
\author[J.L. V\'azquez]{Juan Luis V\'azquez}

\address{\hspace*{-30pt}
{\sc{A. Blanchet:}} \rm CRM, Apartat 50, 08193 Bellaterra, Spain. \Email{blanchet@ceremade.dauphine.fr}\newline
\hspace*{-15pt}{\sc{M. Bonforte:}} \rm Departamento de Matem\'{a}ticas, Univ. Aut\'{o}noma de Madrid, Campus de Cantoblanco, 28049 Madrid, Spain \& Dip. di Matematica, Politecnico di Torino, Corso Duca degli Abruzzi 24, 10129 Torino, Italy.\newline \Email{bonforte@ceremade.dauphine.fr}\newline
\hspace*{-15pt}{\sc{J. Dolbeault:}} \rm CEREMADE, Univ. Paris Dauphine, Pl. de Lattre de Tassigny, 75775 Paris 16, France. \newline \Email{dolbeault@ceremade.dauphine.fr}\newline
\hspace*{-15pt}{\sc{G. Grillo:}} \rm Dip. di Matematica, Politecnico di Torino, Corso Duca degli Abruzzi 24, 10129 Torino, Italy.\newline \Email{gabriele.grillo@polito.it}\newline
\hspace*{-15pt}{\sc{J.L. V\'azquez:}} \rm Departamento de Matem\'{a}ticas, Univ. Aut\'{o}noma de Madrid, Campus de Cantoblanco, 28049 Madrid, Spain. \Email{juanluis.vazquez@uam.es}}

\date{\today}

\begin{abstract}
We consider non-negative solutions of the fast diffusion equation
$u_t=\Delta u^m$ with $m \in (0,1)$, in the Euclidean space
$\RR^d$, $d\ge 3$, and study the asymptotic behavior of a natural
class of solutions, in the limit corresponding to $t\to\infty$ for
$m\ge m_c=(d-2)/d$, or as $t$ approaches the extinction time when
$m<m_c$. For a class   of initial data we prove that the solution
converges with a polynomial rate to a self-similar solution, for
$t$ large enough  if $m\ge m_c$, or close enough to the extinction
time  if $m< m_c$. Such results are new in the range $m\le m_c$
where previous approaches fail. In the range $m_c<m<1$ we improve
on known results.
\end{abstract}

\maketitle
\thispagestyle{empty}

\section{Introduction}\label{Sec:Intro}

We study the Cauchy problem for the fast diffusion equation posed
in the whole Euclidean space, that is, we consider the solutions
$u(\tau,y)$ of
\begin{equation}\label{FDE.Problem}
\left\{\begin{array}{lll}
\partial_\tau u=\Delta u^m\vspace{.3cm}\\
u(0,\cdot)=u_0\,,
\end{array}\right.
\end{equation}
where $m \in (0,1)$ (which means fast diffusion) and
$(\tau,y)\in(0,T)\times\RR^d$ for some $T>0$. We consider
non-negative initial data and solutions. Existence and uniqueness
of weak solutions of this problem with initial data in $\LL1_{\rm
loc}(\RR^d)$ was first proved by M.A. Herrero and M.~Pierre
in~\cite{HP}. In the whole space, the behavior of the solutions is
quite different in the parameter ranges $m_c<m<1$ and $0<m<m_c$,
the critical exponent being defined as
\[
m_c:=\frac{d-2}{d}\;.
\]
Note that $m_c>0$ only if $d\ge 3$, so that the lower range does
not exist for $d=1$, $2$. For $m>m_c$ the mass $\int_{\RR^d}
u(y,t)\dd y$ is preserved in time if the initial datum $u_0$ is
integrable in $\RR^d$. Besides, non-negative solutions are
positive and smooth for all $x\in \RR^d$ and $t>0$. On the
contrary, solutions may extinguish in finite time in the lower
range $m<m_c$, for instance when the initial data is in
$\LL^{p_*}(\RR^d)$ with $p_*=d\,(1-m)/2$: then there exists a time
$T>0$ such that
\begin{equation*}
\lim_{\tau \nearrow T} u(\tau,y)= 0\;.
\end{equation*}
Many computations are however similar in both ranges, from an
algebraic point of view. We refer to the monograph
\cite{Vazquez2006} for a detailed discussion of the existence
theory and references to the subject. The extension
to exponents $m\le 0$ is also treated, and it is natural but it
will not be the focus of this paper.

\medskip

In the last two decades, special attention has been given to the
study of large time asymptotics of these equations, starting with
the pioneering work of A. Friedman and S. Kamin~\cite{fk80} and
completed in \cite{MR1977429}, when $m$ is in the range
$(m_c,\infty)$. In those studies the class of non-negative, finite
mass solutions are considered. Asymptotic stabilization towards
{\sl self-similar asymptotic solutions\/} known as Barenblatt
solutions is shown. For $m_c<m<1$, such solutions take the form:
\begin{equation}\label{baren.form1}
U_{D,T}(\tau,y):=\frac 1{R(\tau)^d} \left(D+ \frac{1-m}{2\,m}
\left|\frac{y}{R(\tau)}\right|^2\right)^{-\frac{1}{1-m}}
\end{equation}
with $R(\tau):=\big[d\,(m-m_c)\,(\tau+
T)\big]^{\frac{1}{d\,(m-m_c)}}$. Here $D, T\ge 0$ are free
parameters. While the second parameter means a time displacement
and does not play much role in the asymptotic behavior, the first
does and can be computed from the mass of the solution. The value
$m_c$ is the critical exponent below which the Barenblatt
solutions cease to exist in this standard form.

Here, we are mainly interested in addressing the question of the
asymptotic behavior of~\eqref{FDE.Problem} when $0<m< m_c$. We
consider a wide class of solutions which vanish in finite time $T$
and describe their behavior as $\tau$ goes to $T$. We point out
that our methods allow to treat simultaneously the ranges $0<m<
m_c$ and $m_c\le m<1$, in which one is interested in the behavior
of the solutions as~$\tau$ goes to infinity. For this purpose, we
extend the Barenblatt solutions to the range $0<m<m_c$ with the
same expression~\eqref{baren.form1}, but a different form for $R$,
that is
\begin{equation*}
R(\tau):=\big[d\,(m_c-m)\,(T-\tau)\big]^{-\frac{1}{d\,(m_c-m)}}\,.
\end{equation*}
The parameter $T$ now denotes the extinction time. Following
\cite{Vazquez2006}, we shall call such solutions the {\sl
pseudo-Barenblatt solutions.\/} Notice that Barenblatt and
pseudo-Barenblatt solutions $U_{D,T}$, with $D$, $T>0$, are such
that $U_{D,T}^p$ is integrable if and only if $p>p_*$ ($p_*$ is
defined above, and $p_*>1$ means $m<m_c$). Consistently with the
above choices, for $m=m_c$, one has to choose
$R(\tau):=e^{\tau+T}$ with free parameter~$T$, see
\cite{Vazquez2006}, in order to obtain pseudo-Barenblatt
solutions; then, $p_*=1$.

The family of Barenblatt (respectively pseudo-Barenblatt)
solutions represents the asymptotic patterns to which many other
solutions converge for large times if $m>m_c$ (respectively as $t$
goes to $T$ if $0<m<m_c$). We are interested in the class of
solutions for which such a convergence takes place and in the
rates of convergence. Both questions strongly depend on $m$. Let
us emphasize for instance that the Barenblatt solution $U_{D,T}$
is integrable in $y$ for $m>m_c$, while the pseudo-Barenblatt
solution corresponding to $m\leq m_c$ is not integrable. Since
much is known in the case $m>m_c$, see for instance
\cite{MR1986060,MR1940370} and
\cite{BV-harnackExtTime,BV,MR1853037,MR1901093,MR1777035,MR2126633,MR1974458,MR1977429}
for more complete results, the main novelty of our paper is
concerned with the lower range $m\le m_c$, which has several
interesting new features. For instance, in the analysis in high
space dimensions, that is $d>4$, another critical exponent
appears,
\begin{equation*}
m_*:=\frac{d-4}{d-2}<m_c\;.
\end{equation*}
A key property of $m_*$ is that the difference of two
pseudo-Barenblatt solutions is integrable for $m\in(m_*,m_c)$,
while it is not integrable for $m\in (0,m_*]$.

The convergence towards Barenblatt and pseudo-Barenblatt solutions
is subtle since the solutions converge to zero everywhere. To
capture the asymptotic profiles, it is therefore convenient to
rescale the solutions and replace the study of intermediate
asymptotics by the study of the convergence to stationary
solutions in {\sl rescaled variables,\/}
\begin{equation}\label{eq:chgvariable}
t:=\log\left(\frac{R(\tau)}{R(0)}\right)\quad\mbox{and}\quad x:=\frac y{R(\tau)}\;,
\end{equation}
with $R$ as above. In these new variables, if $u$ is a solution
to~\eqref{FDE.Problem}, the function
\begin{equation*}
v(t,x):= R(\tau)^{d}\,u(\tau,y)
\end{equation*}
solves a nonlinear {\sl Fokker-Planck type equation,\/}
\begin{equation}\label{FP.Problem}
\left\{
\begin{array}{ll}
\partial_t v(t,x)=\Delta v^m(t,x)+\nabla\cdot(x\,v(t,x)) &
\qquad (t,x) \in(0,+\infty)\times\RR^d\,,\vspace{.3cm}\\
v(0,x)=v_0(x)& \qquad x \in\RR^d\;.
\end{array}
\right.
\end{equation}
The initial data for~\eqref{FDE.Problem} and for the rescaled
equation~\eqref{FP.Problem} are related by
\begin{equation*}
u_0(y)=R(0)^{-d}\,v_0(y/R(0))\;,
\end{equation*}
where $R(0)=\big[d\,|m-m_c|\,T\big]^{\frac{1}{d\,(m-m_c)}}$ only
depends on $T$. In this formulation, the Barenblatt and
pseudo-Barenblatt solutions are transformed into stationary
solutions given by
\begin{equation}\label{baren.form2}
V_D(x):=\left(D+\frac{1-m}{2\,m}\,|x|^2\right)^{-\frac{1}{1-m}}
\end{equation}
where $0<m<1$ and $D> 0$ is a free parameter. With a
straightforward abuse of language, we say that $V_D$ is a {\sl
Barenblatt profile,\/} including the case $m\leq m_c$. The value
$D=0$ can also be admitted as a limit case, but the corresponding
solution is singular at $x=0$. See~\cite{Vazquez2006} for more
details. The parameter $T$ has disappeared from the new problem,
but it enters in the change of variables. Note that in all cases,
$t$ runs from $0$ to infinity in these rescaled variables.

\subsection*{Assumptions and main results}

We can write the assumptions on the {\sl initial conditions\/} in
terms of either $u_0$ or $v_0$. We assume that\smallskip

\noindent (H1) $u_0$ is a non-negative function in $\LL1_{\rm
loc}(\RR^d)$ and that there exist positive constants $T$ and
$D_0>D_1$ such that
\begin{equation*}
U_{D_0,T}(0,y)\le u_0(y)\le U_{D_1,T}(0,y)\quad\forall\;y \in\RR^d\,.
\end{equation*}
\noindent (H2) There exist $D_*\in [D_1,D_0]$ and $f\in\LL1(\RR^d)$ such that
\begin{equation*}
u_0(y)=U_{D_*,T}(0,y)+f(y)\quad\forall\;y \in\RR^d\,.
\end{equation*}
Note that by the Comparison Principle, see Lemma~\ref{Lem:MP}
below, in the case $m<m_c$, (H1) implies that the extinction
occurs at time $T$. When $m>m_*$, (H2) follows from (H1) since the
difference of two Barenblatt solutions is always integrable. For
$m\leq m_*$, (H2) is an additional restriction. We shall assume
throughout this paper that $d\geq 3$ and observe that (H2) has to
be taken into account only if $m_*>0$, that is, $d\geq 5$.

In terms of $v_0$, with $f$ replaced by $R(0)^{-d}\,f(y/R(0))$,
conditions (H1) and (H2) can be rewritten as follows. To avoid
more notations, we keep using $f$ in (H2') although it is not the
same function as in~(H2).\smallskip

\noindent (H1') $v_0$ is a non-negative function in $\LL1_{\rm
loc}(\RR^d)$ and there exist positive constants $D_0> D_1$ such
that
\begin{equation*}
\label{eq:assumptionv}
V_{D_0}(x)\le v_0(x)\le V_{D_1}(x)\quad\forall\;x \in\RR^d\,.
\end{equation*}
(H2') There exist $D_*\in [D_1,D_0]$ and $f\in\LL1(\RR^d)$ such that
\begin{equation*} \label{eq:assumptionmsmallermstarv}
v_0(x)=V_{D_*}(x)+f(x)\quad\forall\;x \in\RR^d\,.
\end{equation*}

If $m\in(m_*,1)$, the map $D\mapsto\int_{\RR^d}(v_0-V_D)\dx$ is
continuous, monotone increasing. Hence we can also define a unique
$D_*\in[D_1,D_0]$ such that
\[
\int_{\RR^d}(v_0-V_{D_*})\dx=0\;.
\]

Before stating any result, one more exponent is needed. We define
$\q(m)$ as the infimum of all positive real numbers $p$ for which
two Barenblatt profiles $V_{D_1}$ and $V_{D_2}$ are such that
$|V_{D_1}-V_{D_2}|$ belongs~to~$\LL^p(\RR^d)$:
\begin{equation*}\label{q_0}
\q(m):=\dfrac{d\,(1-m)}{2\,(2-m)}\;.
\end{equation*}
We see that $\q(m)>1$ if $m\in(0,m_*)$, $\q(m_*)=1$, and $\q(m)<1$ if $m>m_*$.

We can now state the convergence of $v(t)$ towards a unique
Barenblatt profile. For simplicity, we will write $v(t)$ instead
of $x \mapsto v(t,x)$ whenever we want to emphasize the dependence
in $t$.
\begin{theorem}[Convergence to the asymptotic profile]\label{Thm:A1}
Let $d\ge3$, $m\in(0,1)$. Consider the solution $v$
of~\eqref{FP.Problem} with initial data satisfying {\rm
(H1')-(H2')}.
\begin{enumerate}
\item[{\rm (i)}] For any $m>m_*$, there exists a unique
$D_*\in[D_1,D_0]$ such that $\int_{\RR^d}(v(t)-V_{D_*})\dx=0$ for
any $t>0$. Moreover, for any $p\in (\q(m),\infty]$,
$\lim_{t\to\infty}\int_{\RR^d}|v(t)-V_{D_*}|^p\dx=0$. \item[{\rm
(ii)}] For $m\le m_*$, $v(t)-V_{D_*}$ is integrable,
$\int_{\RR^d}(v(t)-V_{D_*})\dx=\int_{\RR^d}f\dx$ and $v(t)$
converges to $V_{D_*}$ in $\LL^p(\RR^d)$ as $t\to\infty$, for any
$p\in (1,\infty]$. \item[{\rm (iii)}] {\rm (Convergence in
Relative Error)} For any $p\in (d/2,\infty]$,
\begin{equation*}\label{CRE}
\lim_{t\to\infty}\left\|\,{v(t)}/\,{V_{D_*}}-1\,\right\|_{p}=\;0\;.
\end{equation*}
\end{enumerate}
\end{theorem}
\noindent In case $m>m_*$, the value of $D_*$ can be computed at
$t=0$ as a consequence of the mass balance law
$\int_{\RR^d}(v_0-V_{D_*})\dx=0$, and then the conservation result
holds for all $t>0$, see Proposition~\ref{prop:relconsmass} below.
On the other hand, in the case $m\le m_*$ the mass balance does
not make sense, but $D_*$ is determined by Assumption (H2). In
this case, the presence of a perturbation $f$ of $V_{D_*}$ with
nonzero mass, does not affect the asymptotic behavior of the
solution at first order.

In a recent paper \cite{Daskalopoulos-Sesum2006}, P.~Daskalopoulos
and N.~Sesum prove some of the results of Theorem~\ref{Thm:A1}
under similar hypotheses (see \cite[Theorem
1.4]{Daskalopoulos-Sesum2006}). Actually they only prove the
$\LL^\infty$ convergence in case (ii) and the
$\LL1\cap\LL^\infty$ convergence in case (i). Our proof was
obtained independently and announced in~\cite{BBDGV-CRAS}. It is
based on entropy estimates and paves the way to the sharper
results on convergence with rates, which are the main purpose of
the present paper. Assertion (iii) says that the convergence of
(i)-(ii) can be improved into a convergence in relative error, in
the sense of~\cite{MR1977429}. Such a strong convergence may look
surprising at first sight, but it is a consequence of Assumption
(H1'): the tails of $v_0$ and $V_{D_*}$ have the same behavior as
$|x|\to\infty$.

We can now state our main asymptotic result, on rates of
convergence. To state this second result, we need yet another
exponent,
\begin{equation*}\label{q_*}
q_*:= \dfrac{2\,d\,(1-m)}{2\,(2-m)+d\,(1-m)}
\end{equation*}
and note that $q_*>1$ if $m<m_*$, $q_*=1$ if $m=m_*$, and $q_*<1$
if $m>m_*$. For any $q>q_*$, the function $V_{D_*}$ is in
$\LL^{(2-m)q/(2-q)}(\RR^d)$, which allows us to use convenient
H\"older interpolation inequalities. We define the $C^j$ semi-norm
by
\begin{equation*}
\|f\|_{C^j(\RR^d)}:=\displaystyle \max_{|\eta|=j}\; \displaystyle
\sup_{x\in\RR^d}\big|\partial^\eta f(x)\,\big|\,,
\end{equation*}
where the standard multi-index notation is used:
$|\eta|=\,\eta_1+\,\ldots\,+\eta_d\;$ is the length of the
multi-index $\eta=(\eta_1,\,\ldots\,,\eta_d)\in\mathbb{Z}^d$. The
last ingredient is a Hardy-Poincar\'e constant, which is defined
as follows. For any $m\in(0,m_*)\cup(m_*,1)$, let
\begin{equation}\label{Defn-Lambda}
\gopt\;:=\;m\;\inf_h\;\frac{\int_{\RR^d}|\nabla
h|^2\,V_{D_*}\dx}{\int_{\RR^d}|h-\bar h|^2\,V_{D_*}^{2-m}\dx}\;,
\end{equation}
where the infimum is taken over the set of smooth functions $h$
such that ${\rm supp}(h)\subset\RR^d\setminus\{0\}$ and $\bar h=0$
if $m<m_*$, while for $m>m_*$, $\bar
h:={\int_{\RR^d}h\,V_{D_*}^{2-m}\dx}/{\int_{\RR^d}V_{D_*}^{2-m}\dx}$,
{\sl \cf\/}Theorem~\ref{Teor.Spect.} in the Appendix for more
details. We shall prove that $\gopt$ is positive and independent
of $D_*$. In the next result, the time decay rate is formulated in
terms of the spectral gap $\gopt$. Analyzing the relationship
between the optimal constant $\mathcal C_{m,d}=m/\gopt$ in the
corresponding functional inequality and the asymptotic rates of
the fast diffusion equation is the leitmotiv of this paper. The
case $m=m_*$ has to be excluded for reasons which are deeply
related to Hardy's inequality, see~\cite{BBDGV-CRAS}.
\begin{theorem}[Convergence with rate]\label{Thm:A2}
Under the assumptions of Theorem~{\rm~\ref{Thm:A1}}, with $\gopt$
given by~\eqref{Defn-Lambda}, if $m\ne m_*$, there exists $t_0\geq
0$ such that the following properties hold:
\begin{enumerate}
\item[{\rm (i)}] For any $q\in(q_*,\infty]$, there exists a
positive constant $C_q$ such that
\begin{equation*}
\|v(t)-V_{D_*}\|_q\le C_q\;\ee^{-\gopt \,t}\quad\forall\;t\geq t_0\;.
\end{equation*}
\item[{\rm (ii)}] For any $ \vartheta\in[0,(2-m)/(1-m))$, there
exists a positive constant $K_\vartheta$ such that
\begin{equation*}\label{Moments}
\big\|\,|x|^\vartheta(v(t)-V_{D_*})\big\|_{2}\le
K_\vartheta\;\ee^{-\gopt\,t}\quad\forall\;t\geq t_0\;.
\end{equation*}
\item[{\rm (iii)}] For any $j\in\mathbb{N}$, there exists a positive constant $H_j$ such that
\begin{equation*}\label{Estim.C^j.Intro}
\|v(t)-V_{D_*}\|_{C^j(\RR^d)}\le H_j\,\ee^{-\frac{2\,\gopt}{d+2(j+1)}\,t}\quad\forall\;t\geq t_0\;.
\end{equation*}
\end{enumerate}
\end{theorem}
The constants $C_q$, $K_\vartheta$ and $H_j$ depend on $t_0$, $m$,
$d$, $v_0$, $D_0$, $D_1$, and $q$, $\vartheta$ and $j$; $t_0$ also
depends on $D_0$ and $D_1$. It is remarkable that the decay rate
of the nonlinear problem is given exactly by $\gopt$ (see
Section~\ref{sec.proof.CRE-exp}). Using~\eqref{eq:chgvariable},
the results of Theorem~\ref{Thm:A2} for the solution $v(t)$
of~\eqref{FP.Problem} can be translated into results for the
solution $u(\tau)$ of~\eqref{FDE.Problem} as follows.
\begin{corollary}[Intermediate asymptotics]\label{Cor:A2}
Let $d\ge3$, $m\in(0,1)$, $m\ne m_*$. Consider a solution~$u$
of~\eqref{FDE.Problem} with initial data satisfying {\rm
(H1)-(H2)}. For $\tau$ large enough, for any $q \in (q_*,\infty]$,
there exists a positive constant $C$ such that
\begin{equation*}
\|u(\tau)-U_{D_*}(\tau)\|_q\le C\,R(\tau)^{-\alpha}\;,
\end{equation*}
where $\alpha=\gopt+d\,(q-1)/q$ with $\gopt$ given
by~\eqref{Defn-Lambda}, and {\rm large} means $T-\tau >0$, small,
if $m<m_c$, and $\tau\to\infty$ if $m\ge m_c$.
\end{corollary}
We also obtain a convergence result in relative error. For any $p\in(d/2,\infty]$, define
\begin{equation*}
\lambda(p):=\frac{(2\,p-d)\,(1-m)}{p\,(d+2)\,(2-m)}\,\gopt\;.
\end{equation*}
\begin{theorem}[Exponential Decay of Relative Error]\label{thm:CRE-exp}
Under the assumptions of Theorem~{\rm \ref{Thm:A2}}, if $m\ne
m_*$, for any $p\in(d/2,\infty]$, there exists a positive constant
$\mathcal{C}$ and $\lambda\in (0,\lambda(p))$ such that
\begin{equation*}
\big\|\,{v(t)}/\,{V_{D_*}}-1\,\big\|_p\le \mathcal{C}\,\ee^{-\lambda\,t}\quad\forall\;t\geq 0\;.
\end{equation*}
\end{theorem}
Let us list a few observations on the above results.

\noindent (a) In the range $m_c<m<1$, convergence with rates has
been obtained under various assumptions,
\cf~\cite{MR1940370,MR1853037,MR1777035} (optimal rates) for
$m\in[m_1,1)$, $m_1=(d-1)/d$, and \cite{MR1901093,MR1986060} for
$m\in(m_c,m_1)$, which are weaker than the ones of
Theorem~\ref{Thm:A2}. See \cite{MR2126633} for the detailed
analysis of the spectrum of the linearized operator in the range
$m>m_c$. A stronger convergence has also been proved in the sense
of relative error under very mild assumptions,
\cf~\cite{MR1977429}.

\noindent (b) In the range $m_c<m<1$, Assumption (H1) is less
restrictive than one could think. By the global Harnack principle
of \cite{BV}, see Theorem~\ref{thm.ghp} below, any solution with
non-negative initial data $u_0\in \LL1_{\rm loc}(\RR^d)$ that
decays at infinity like $u_0(y)= O(|y|^{2/(1-m)})$, is indeed
trapped for all $t>0$ between two Barenblatt solutions if
$m_c<m<1$. The restrictions on the class of initial data are
therefore not so essential as far as the asymptotic behavior is
concerned, and can therefore be relaxed.

\noindent (c) In the range $0<m\le m_c$, the pseudo-Barenblatt
solutions are not integrable. For $m<m_c$ many solutions vanish in
finite time and have various asymptotic behaviors depending on
the initial data. Solutions with bounded and integrable initial
data are described by self-similar solutions with so-called
anomalous exponents, see~\cite{King,PelZhang}
and~\cite[Chapter~7]{Vazquez2006}. Even for solutions with initial
data not so far from a pseudo-Barenblatt solution, the asymptotic
behavior may significantly differ from the behavior of a
pseudo-Barenblatt solution: in~\cite[Theorem
1.4]{Daskalopoulos-Sesum2006} a solution of~\eqref{FDE.Problem} is
found, which is such that
$\lim_{|x|\to\infty}u_0(x)/U_{D,T}(0,x)=1$ and which, after
rescaling, does not converge to $V_D$ as $t\to\infty$, that is for
$\tau\to T$. Assumption (H1)-(H2) are therefore more restrictive
than for $m>m_c$.

\noindent (d) Proofs are constructive and the values of the
various constants are explicit although not so easy to write. The
interested reader will be able to recover their expressions by
carefully reading the proofs, where all details are given. See
Appendix A for more details on the constant $\gopt$ that controls
the rate of convergence. The rate given by this exponent is sharp
in the linear case, and a deeper analysis should prove that it is
sharp also in the nonlinear case. Obtaining the optimal value of
$\gopt$ is still an open question in the range $m\in(m_*,m_c]$.
Our method gives convergence with rates even in the limit case
$m=m_c$, which is new.

\subsection*{Further comments}

After stating our main results, let us come back a little bit on
the motivations of this paper, on the main tools and the
originality of our results with respect to the existing
literature.

During the last few years, asymptotic rates of convergence for the
solutions of nonlinear diffusion equations have attracted lots of
attention, usually in connection with time-dependent scalings and
entropy methods. This has been first done in the range of
exponents corresponding to the porous medium equation, with
$1<m<2$, and in the range where standard Gagliardo-Nirenberg
inequalities apply, $m_1\le m<1$,
see~\cite{MR1940370,MR1853037,MR1777035}. The class of
non-negative, finite mass solutions has to be narrowed to the
smaller set of functions with finite free energy, or to be
precise, with finite entropy and finite potential energy. In the
rescaled variables, asymptotic stabilization to the Barenblatt
profiles holds at an exponential rate, while in the original time
variable $\tau$, the convergence of the difference with the
Barenblatt solutions holds at a power-law rate, which is shown to
be optimal.

The next question was to understand what happens for $m<m_1$. After the linearized analysis of~\cite{MR1901093}, the proof of convergence with rates was done in~\cite{MR1986060} in the range $ m_c<m<m_1$ for which global existence of finite mass solutions still holds. The basin of attraction is narrowed to the class of solutions with finite relative entropy with respect to some Barenblatt solution.

A dramatic change occurs for $m<m_c$, since a large class of solutions vanish in finite time. As a consequence, mass is not conserved, and a key estimate for higher values of $m$ is lost. It is however natural to investigate the basin of attraction of the pseudo-Barenblatt solutions for $m\le m_c$ using relative entropy techniques and to study the convergence rates. This can be done in a weighted space using functional inequalities, which can still be related to some spectral properties of a differential operator obtained by an appropriate linearization.

\smallskip The {\sl generalized entropy functional,\/} or {\sl free energy functional,\/} is defined as
\begin{equation*}
\fe{v}:=\int_{\RR^d}\Big[\varphi(v)+\frac 12\,|x|^2\,v\Big]\dx\quad \mbox{where}\quad \varphi(v):=\frac{v^m}{m-1}\;.
\end{equation*}
It is then observed that the free energy of the Barenblatt profiles, \cf~\cite{MR1777035,MR1940370}, becomes infinite if $m\leq m_0, $ where $m_0:=d/(d+2)\in(m_c,m_1)$. In order to avoid this difficulty, it is convenient to work with the {\sl relative entropy\/} of $v$ with respect to $V_D$ defined as follows:
\begin{equation*}
\fe{v|V_D}:=\int_{\RR^d}\Big[\varphi(v)-\varphi(V_D)- \varphi'(V_D)\,(v-V_D)\Big]\dx\;.
\end{equation*}
The {\sl relative entropy\/} is the key tool of our analysis. It is such that $\fe{v|V_D}:=\fe{v}-\fe{V_D}$ if $m\in(m_0,1)$ and $\int_{\RR^d}v\dx= \int_{\RR^d}V_D\dx$, that is for $D=D_*$, where $D_*$ is as in Theorem~\ref{Thm:A1}. The functional $\fe{v|V_{D_*}}$ can also be defined for $m\leq m_0$. By homogeneity of $\varphi$, we can indeed rewrite it as
\[
\fe{v|V_{D_*}}:= \int_{\RR^d}\Big[\varphi(w)-\varphi(1)- \varphi'(1)\,(w-1)\Big]\,V_{D_*}^m\dx\quad\mbox{with}\quad w=\frac v{V_{D_*}}\;.
\]
This makes clear why it is well defined at least for $w$ close enough to $1$ as $|x|\to\infty$. The functional $v\mapsto\fe{v|V_{D_*}}$ is convex and achieves its minimum, $0$, for $v=V_{D_*}$. If $v$ is a solution of \eqref{FP.Problem}, the {\sl entropy production term\/} takes the form
\begin{equation*}
- \frac d{dt} \fe{v(t)|V_{D_*}}=\I{v(t)}{V_{D_*}}\;,
\end{equation*}
where the functional
$$
v\mapsto\I{v}{V_D}:=\int_{\RR^d}v\,\Big|\nabla\varphi'(v)-\nabla\varphi'(V_D)\Big|^2\dx
$$
will be called {\sl the relative Fisher information}. See Proposition~\ref{Prop:EntrProd} for more details. For any $m\in [m_1,1)$, $\fe{v|V_{D_*}} \le\frac 12\,\I{v}{V_{D_*}}$ holds for any smooth function $v$ and the inequality is nothing else than the optimal Gagliardo-Nirenberg inequality, for which equality is achieved precisely by the Barenblatt profiles, see~\cite{MR1940370}. In such a case,
\[
\fe{v(t)|V_{D_*}}\leq \fe{v_0|V_{D_*}}\,e^{-\,2\,t}\quad\forall\; t\geq 0\;.
\]
The limit case $m=m_1$ corresponds to the critical Sobolev
inequality whose optimal form was established by T. Aubin and G.
Talenti in~\cite{MR0448404,MR0463908}, while in the limit $m\to 1$
one recovers Gross' logarithmic Sobolev inequality,
see~\cite{Gr-LSI,MR1940370}. For $m\in[m_1,1)$, F. Otto
in~\cite{MR1842429} noticed that~\eqref{FP.Problem} can be
interpreted as the gradient flow of the free energy with respect
to the Wasserstein distance. The exponent $m=m_1$ is the limit
case for which the {\sl displacement convexity\/} property holds
true.

Pushing the method to the case $0<m<m_1$ requires the use of the
relative entropy in place of the free energy. The method applies
only to a class of initial data which have a finite relative
entropy with respect to some Barenblatt profile $V_{D_*}$ and
satisfy convenient bounds. Mass can be finite in the case
$m\in(m_c,m_1)$, which was the framework of some earlier studies,
see~\cite{MR1901093,MR1986060}, or infinite if $m\in (0,m_c)$. Two
Barenblatt profiles $V_{D_0}$ and $V_{D_1}$ have finite relative
entropy, \ie\ $\fe{V_{D_1}|V_{D_0}}<\infty$ if and only if either
$d\leq 4$, or $d\geq 5$ and $m>m_*$, $m_*=(d-4)/(d-2)$. Hence, for
$d\geq 5$, $m=m_*$ is a threshold not only for defining the
relative mass of two pseudo-Barenblatt solutions, but also for
defining their relative entropies or for the integrability of
$V_{D_*}^{2-m}$. Note that $m_*<m_c$ for all~$d\geq 5$. The proof
of Theorem~\ref{Thm:A2} amounts to prove that the relative entropy
$\fe{v|V_{D_*}}$ decays in time and converges to $0$ at an
exponential rate when $t\to\infty$. For $m>\min\{0,m_*\}$,
$\fe{v|V_{D_*}}$ is well defined under condition (H1'). For
$m<m_*$, an additional restriction is required, which is precisely
the purpose of (H2').

Our approach of course covers the case $m\ge m_c$ and we recover
some of the results found in~\cite{MR1901093,MR1986060}. Some of our
results can also be extended to the range $m<0$, but additional
technical complications arise, which are still to be studied. In
this paper, we leave apart several interesting questions, like the
precise study of the case of $m=m_*$ or the equation $u_t=\Delta
\log u$ in dimension $d\ge 2$, see {\sl e.g.\/}
\cite{MR1679782,Daskalopoulos-Sesum,MR1606339,MR1357953}, which is
the natural limiting equation to study in the limit $m\to 0$. Also
see \cite{MR2162628,MR2194833,MR2145602} for results which seem
closely related to ours, and \cite{MR1857048} in the case
$m=(d-2)/(d+2)$. In particular we do not use the Bakry-Emery method
introduced in \cite{MR889476}, on which the results of
\cite{MR1777035,MR1853037,MR1974458,MR1901093,MR1986060} are based.
We prove a conservation of relative mass, which allows us to remove
the limitation $m>m_c$. Neither mass transportation techniques nor
Wasserstein distance are needed, although the approach of
Section~\ref{Sec:A2} is not unrelated, see
\cite{bobkov-gotze,BR1,BR2,Miclo}.

 The paper is organized as follows. In Section \ref{Sec:Basic}, we
extend the property of mass conservation, which holds only for
$m>m_c$, to a property of conservation of relative mass, see
Proposition~\ref{prop:relconsmass}. This selects a unique
Barenblatt profile, which governs the asymptotic behavior. We also
establish regularity properties of the solutions. From there on we
work with the quotient of the solution by the Barenblatt profile.
In Section~\ref{Sec:ConvergenceWithoutRate}, we prove
Theorem~\ref{Thm:A1} and establish several properties which are
used in the sequel. Sections~\ref{Sec:SpectralGap}
and~\ref{sec:cvrgcewithoutrate} are respectively devoted to
introducting a suitable linearization and to the derivation of
entropy - entropy production estimates in the nonlinear case, from
the corresponding spectral gap estimates for the linearized
problem. The proof of Theorem~\ref{Thm:A2} is given in
Section~\ref{Sec:WithRate}.

Appendix A is devoted to the proofs of spectral gap estimates,
that is, of weighted Poincar\'e-Hardy inequalities, which have
already been studied in~\cite{BBDGV-CRAS}, and in
\cite{MR1901093,MR2126633} in the special case \hbox{$m>m_c$}. We
consider the family of weights of the form $V_D$ or $V_D^{2-m}$,
$D>0$, that are obtained from the linearization of the relative
entropy. In the limit $D\to 0$, they yield the case corresponding
to the weighted $\LL2$ norm of the Caffarelli-Kohn-Nirenberg
inequalities, \cf~\cite{CKN, CW}.

A final section, Appendix B, explains how to extend the results of
this paper to the fast diffusion with exponents $m\le 0$. Note
that the equation needs to be properly modified. The conclusion is
that the results still hold and the proofs need only minor
modifications that are indicated.

\section{Basic Estimates}\label{Sec:Basic}

We establish in this section the main properties of the solutions that will be used in the sequel.

\subsection{$\LL1$-contraction and Maximum Principle}

\begin{lemma}[$\LL1$-contraction]\label{Lem:contraction} For any two non-negative solutions $u_1$ and $u_2$ of~\eqref{FDE.Problem} defined on a time interval $[0,T)$, with initial data in $\LL1_{{\rm loc}}(\RR^d)$, and any two times $t_1$ and $t_2$ such that $0\le t_1\le t_2< T$, we have
\begin{equation*}
\int_{\RR^d}|u_1(t_2)-u_2(t_2)|\dx\le \int_{\RR^d}|u_1(t_1)-u_2(t_1)|\dx\;.
\end{equation*} \end{lemma}
The above result is well-known to be true for solutions with $\LL1$ data, \cf~\cite[Proposition 9.1]{VazBook}, even in the stronger form
$$
\int_{\RR^d}\big[u_1(t_2)-u_2(t_2)\big]_+\dx\le \int_{\RR^d}\big[u_1(t_1)-u_2(t_1)\big]_+\dx\;,
$$
where $[u]_+$ denotes the positive part of $u$. The result for
data in $\LL1_{{\rm loc}}(\RR^d)$ is obtained by approximation,
using the uniqueness of solutions to the Cauchy Problem, which has
been established in~\cite{HP}. Note that when the right-hand side
is infinite the result applies but there is nothing to prove. As a
consequence, we also obtain the following.
\begin{lemma}[Comparison Principle]\label{Lem:MP} For any two non-negative solutions $u_1$ and $u_2$ of~\eqref{FDE.Problem} on $[0,T)$, $T>0,$ with initial data satisfying $u_{01}\le u_{02}$ a.e, $u_{02}\in\LL1_{{\rm loc}}(\RR^d)$, then we have $u_1(t)\le u_2(t)$ for almost every $t\in [0,T)$. \end{lemma}
We will see below that under Assumption (H1)-(H2) the solutions are smooth functions, hence the comparison in the previous result holds everywhere in $[0,T)\times \RR^d$.

\subsection{Conservation of relative mass}\label{Subsec:RelMass}

Mass conservation is used in the range $m>m_c$ to determine the
parameter~$D$ which characterizes the Barenblatt profile $V_D$. In
the range $m\leq m_c$, we can still prove that
$\int_{\RR^d}(v(t)-V_{D})\dx$ is conserved for any $t>0$, even if
$V_D\not\in \LL1(\RR^d)$.
\begin{proposition}\label{prop:relconsmass}
Let $m\in (0,1)$. Consider a solution $u$ of~\eqref{FDE.Problem}
with initial data $u_0$ satisfying {\rm (H1)-(H2)}. If for some
$D>0$, $\int_{\RR^d}(u_0-U_{D,T}(0,\cdot))\dd x$ is finite, then
\begin{equation*}\label{integral.eq}
\int_{\RR^d}\big[u(\tau,x)-U_{D,T}(\tau,x)\big]\dd
x=\int_{\RR^d}\big[u_0(x)-U_{D,T}(0,x)\big]\dd x\quad
\forall\;\tau\in(0,T)\;.
\end{equation*} \end{proposition}

\proof In the range $m>m_c$, $u_0$ and $U_{D,T}(0,\cdot)$ are
integrable and mass conservation of the solutions
of~\eqref{FDE.Problem} is well-known.

Assume next that $m<m_c$ and let~$\chi$ be a $C2$ function
on~$\RR^+$ such that $\chi\equiv 1$ on $[0,1]$, $\chi\equiv 0$ on
$[2,\infty)$, and $0\leq\chi\leq 1$ on $[1,2]$. For any
$\lambda>0$, take $\phi_\lambda(x):=\chi(|x|/\lambda)$ as a test
function and denote by $B_\lambda$ the ball $B(0,\lambda)$. Then,
\begin{align*}
\left|\dfrac{\rd}{\rd\tau}\int_{\RR^d}\left[u(\tau)-U_{D,T}(\tau)\right]\,\phi_\lambda\dd y\right| &=\left|\int_{\RR^d}\left[u^m(\tau)-U_{D,T}^m(\tau)\right]\,\Delta\phi_\lambda\dd y\right|\\
&=\left|\int_{B_{2\lambda}\setminus B_\lambda}\left[u^m(\tau)-U_{D,T}^m(\tau)\right]\,\Delta\phi_\lambda\dd y\right|\\
&\quad\le C\, \int_{B_{2\lambda}\setminus B_\lambda}\left|u(\tau)-U_{D,T}(\tau)\right|\,U_{D_0,T}^{m-1}(\tau)\,\left|\Delta\phi_\lambda\right|\dd y\;,
\end{align*}
where $C$ is a numerical constant depending on $D$, $D_0$, $D_1$.
As $\lambda\to\infty$, we observe that in $B_{2\lambda}\setminus
B_\lambda$, $U_{D_0,T}^{m-1}$ and $|\Delta\phi_\lambda|$ behave
like $\lambda2$ and $\lambda^{-2}$, so that
$U_{D_0,T}^{m-1}(\tau)\,\Delta\phi$ is bounded uniformly with
respect to $\lambda$. The right hand side is therefore bounded by
$\int_{B_{2\lambda}\setminus
B_\lambda}\left|u(\tau)-U_{D,T}(\tau)\right|\dd y$. For any
$\tau_1$, $\tau_2\in [0,T)$, we write
\begin{equation*}\label{eq2.2}
\left| \int_{\RR^d}\left[u(\tau_2)-U_{D,T}(\tau_2)\right]\,
\phi_\lambda\dd y - \int_{\RR^d}\left[u(\tau_1)-U_{D,T}(\tau_1)\right]\,
\phi_\lambda\dd y \right|\le\int_{\tau_1}^{\tau_2}\kern-4pt\int_{B_{2\lambda}\setminus B_\lambda}\kern-24pt|u-U_{D,T}|\dd y\dd\tau\;.
\end{equation*}
By the $\LL1$-contraction principle, see
Lemma~\ref{Lem:contraction}, we also know that
$|u(\tau,y)-U_{D,T}(\tau,y)|$ is integrable in $y$, uniformly for
all $\tau>0$. The integrability condition implies that the
right-hand side goes to zero in the limit $\lambda\to +\infty$.

\noindent The case $m=m_c$ is similar, except that there is no
extinction time. \qed

In the rescaled variables given by~\eqref{eq:chgvariable},
relative mass is also conserved. Consider a solution $v$
of~\eqref{FP.Problem} with initial data $v_0$ satisfying {\rm
(H1')-(H2')}. If for some $D>0$,
$\int_{\RR^d}(v_0(x)-V_{D}(t,x))\dd x$ is finite, then
\[
\int_{\RR^d}\big[v(t,x)-V_{D}(t,x)\big]\dd
x=\int_{\RR^d}\big[v_0(x)-V_{D}(t,x)\big]\dd x\quad \forall\;
t>0\;.
\]

Whenever $m \le m_c$, we recall that pseudo-Barenblatt solutions
have infinite mass, that is \hbox{$\int_{\RR^d}V_{D_*}\dd
y=\infty$}, but we observe that the difference of two
pseudo-Barenblatt solutions is integrable if $m>m_*$. In such a
case, the parameter $D\in[D_1,D_0]$ in
Proposition~\ref{prop:relconsmass} can be arbitrary. In the proof,
we can moreover estimate $\int_{B_{2\lambda}\setminus
B_\lambda}|u-U_{D,T}|\dd y$ by
\[
\max\left\{\int_{B_{2\lambda/R(t)}\setminus B_{\lambda/R(t)}}|V_{D_0}-V_D|\dd y\;,\quad\int_{B_{2\lambda/R(t)}\setminus B_{\lambda/R(t)}}|V_{D_1}-V_D|\dd y\right\}\;,
\]
which converges to $0$ as $\lambda\to\infty$. As already quoted in the introduction, the map $D\mapsto\int_{\RR^d}(v_0-V_D)\dx$ is continuous, monotone increasing, and we can define a unique $D_*\in[D_1,D_0]$ such that $\int_{\RR^d}(v_0-V_{D_*})\dx=0$. Under the assumptions of Proposition~\ref{prop:relconsmass},
\[
\int_{\RR^d} [v(t,x)-V_{D_*}(x)]\dd x = 0\quad\forall\; t>0\;.
\]
This fact is used in the statement and proof of Theorem \ref{Thm:A1} for $m>m_*$. On the contrary, if $m\in(0,m_*]$, integrals are infinite unless $D=D_*$ in Proposition~\ref{prop:relconsmass}, and then, with the notations of Assumption (H2'),
\[
\int_{\RR^d} [v(t,x)-V_{D_*}(x)]\dd x = \int_{\RR^d} f\dd x\quad\forall\; t>0\;.
\]
In both cases, that is for any $m\in(0,1)$, we shall summarize the fact that $\frac d{dt}\int_{\RR^d} [v-V_{D_*}]\dd x=0$ by saying that the {\sl relative mass\/} is conserved.

\subsection{Passing to the quotient}

Consider a solution $v$ of~\eqref{FP.Problem}. As in~\cite{MR1977429,BV,BV-harnackExtTime}, define
\begin{equation}\label{def.w}
w(t,x):=\frac{v(t,x)}{V_{D_*}(x)}\quad\forall\;(t,x)\in(0,\infty)\times\RR^d\,.
\end{equation}
Next, we rewrite Problem~\eqref{FP.Problem} in terms of $w$:
\begin{equation}\label{FP.Problem.Quotients}
\left\{ \begin{array}{lll} \displaystyle
w_t=\frac{1}{V_{D_*}}\,\nabla\cdot\!\left[w\,V_{D_*}\nabla\!\left(\frac{m}{m-1}\,(w^{m-1}-1) \,V_{D_*}^{m-1}\right)\right]\; &{\rm in}\;(0,+\infty)\times\RR^d\,,\vspace{.3cm}\\
\displaystyle w(0,\cdot)=w_0:=\frac{v_0}{V_{D_*}}\;&{\rm in}\;\RR^d\,.
\end{array} \right. \end{equation}
Define
\begin{equation*}
W_0:= \inf_{x\in\RR^d}\frac{V_{D_0}}{V_{D_*}}\le \sup_{x\in\RR^d}\frac{V_{D_1}}{V_{D_*}}:=W_1\;.
\end{equation*}
A straightforward calculation gives
\begin{equation*}
W_0=\left(\frac{{D_*}}{D_0}\right)^{\frac{1}{1-m}} <\quad 1\quad<\quad \left(\frac{{D_*}}{D_1}\right)^{\frac{1}{1-m}}:=W_1\;.
\end{equation*}
In terms of $w_0$, assumptions (H1') and (H2') can be rewritten as follows:

\smallskip\noindent (H1'') $w_0$ is a non-negative function in $\LL1_{\rm loc}(\RR^d)$ and there exist positive constants $D_0>D_1$ such that
\begin{equation*}
0<W_0\le\frac{V_{D_0}(x)}{V_{D_*}(x)}\le w(x)\le \frac{V_{D_1}(x)}{V_{D_*}(x)}\le W_1<+\infty\quad\forall\;x \in\RR^d\,.
\end{equation*}
(H2'') There exists $f\in \LL1(\RR^d)$ such that
\begin{equation*}
w(x)=1+\frac{f(x)}{V_{D_*}(x)}\quad\forall\;x \in\RR^d\,.
\end{equation*}
As a consequence of the Comparison Principle, see Lemma~\ref{Lem:MP}, (H1'') is satisfied by a solution~$w$ of~\eqref{FP.Problem.Quotients} if it is satisfied by $w_0$.

\subsection{Regularity estimates and Harnack principle}

We start by briefly recalling some well-known results for solutions of fast diffusion equations. A basic regularity result is due to DiBenedetto et al.\,, see {\cite[p. 270]{MR1158938}}, and concerns local space-time H\"older regularity for Problem~\eqref{FDE.Problem}, with some H\"older exponent $\alpha\in(0,1)$; it holds for locally bounded initial data, possibly with sign changes. In the present situation of locally bounded and positive initial data, it is known that the solutions are $C^\infty$ as long as they do not vanish identically because we avoid any degeneracy or singularity and the standard parabolic theory applies. In the sequel we are in particular interested in some sort of uniform $C1$ regularity under the assumption (H1''). We find that it is preferable to work with the function $w$ introduced in~\eqref{def.w} since it is uniformly bounded from above, and from below away from zero.
 --
\begin{theorem}[Uniform $C^k$ regularity]\label{lem:holdereg} Let $m \in (0,1)$ and $w \in \LL^{\infty}_{{\rm loc}}((0,T) \times \RR^d)$ be a solution of~\eqref{FP.Problem.Quotients}. Then for any $k\in\N$, for any $t_0\in (0,T)$,
\begin{equation*}
\sup_{t\geq t_0}\|w(t)\|_{C^k(\RR^d)}<+\infty\;.
\end{equation*} \end{theorem}
\proof Take $t \ge t_0>0$. For a given $\lambda>0$, the equation
for~$w$ is uniformly parabolic in $B_\lambda$, so we conclude that
the regularity estimates hold on $B_{\lambda/2}$ for any $t\ge
t_0$. Let $v$ be the solution of (\ref{FP.Problem}) corresponding
to $w$. In order to cover the large values of $|x|$, we consider
the scaling
\begin{equation*}
v_\lambda(x,t)=\lambda^{2/(1-m)}\,v(\lambda\,x,t)
\end{equation*}
with $\lambda\to\infty$. Then $v_\lambda$ is again a solution
of~\eqref{FP.Problem}, but the region $\Omega_\lambda=\{x\in\RR^d
: \lambda\le |x|\le 2\,\lambda\}$ gets mapped into the region
$\Omega_1=\{x\in\RR^d : 1\le |x|\le 2\}$, for all $\lambda$. Note
also that this scaling transforms the Barenblatt profiles
according to
$$
(V_{D})_{\lambda}=V_{D/\lambda2}\;,
$$
so that, on $\Omega_1$, $(V_{D})_{\lambda}$ is uniformly bounded
from above and from below in $\Omega_1$ as $\lambda\to\infty$, and
converges to $V_0$.

Next, we pass to the functions
$w_\lambda(x,t)=v_\lambda(x,t)/(V_{D_*})_\lambda(x)$ and observe
that in $\Omega_1$, $(w_\lambda)_{\lambda\ge1}$ is uniformly
bounded from above, and from below away from zero. Since
$w_\lambda$ satisfies \eqref{FP.Problem.Quotients} with $V_{D_*}$
replaced by $(V_{D_*})_\lambda$, we conclude as in part~(i) that
$$
\|w_\lambda(t)\|_{C^k(\Omega_1
)}\le C_k
$$
uniformly in $t\ge t_0$ and $\lambda\ge 1$. Undoing the scaling, we find a constant $C$ independent of $\lambda$ such that, for any $t\ge t_0$,
$$
\frac{|\nabla w(t,\lambda\,x)|}{w(t,\lambda\,x)}\le \frac C{\lambda}\;.
$$
We conclude that the result holds for $k=1$. The same argument applies for $k>1$. \qed

In the range $m\in(m_c,1)$, the regularization effects of the fast diffusion equation allow to get rid of assumption (H1). See \cite{BV} for details expressions of the constants entering the statement.
\begin{theorem}[Global Harnack principle, {\cite[Theorem 1.2]{BV}}] \label{thm.ghp} Let $u_0$ be a non-negative function in $\LL1(\RR^d)$ solution of~\eqref{FDE.Problem}. Assume that for some $R>0$, $\sup_{|y|>R}u_0(y)\,|y|^{2/(1-m)}$ is finite. Then, for any $\varepsilon>0$, there exist positive constants $T_0$, $T_1$, $D_0$ and $D_1$, such that
\begin{equation*}
U_{D_0, T_0}(\tau,y)\le u(\tau,y) \le U_{D_1,T_1}(\tau,y)\quad\forall\;(\tau,y)\in (\varepsilon,\infty)\times\RR^d\,.
\end{equation*} \end{theorem}

\subsection{Relative entropy}
\label{rel.entr}

In terms of $w$, we define the {\sl relative entropy\/}
\begin{equation}\label{Entropy.Quotients}
\fw w:=\frac1{1-m}\int_{\RR^d}\big[\,(w-1)-\frac 1m\,(w^m-1)\big]\,V_{D_*}^m\dx
\end{equation}
and the {\sl relative Fisher information\/}
\begin{equation}\label{Fisher.Quotients}
\iw w:=\frac m{(m-1)2}\int_{\RR^d} \big|\,\nabla\left[\left(w^{m-1}-1 \right)V_{D_*}^{m-1}\right]\big|^{\kern 1pt 2}\,w\,V_{D_*}\dx\;.
\end{equation}
These definitions are consistent with the ones given in the introduction, in the sense that, for $w=v/\,V_{D_*}$,
\begin{equation*}
\fw w=\frac 1m\,\fe{v|V_D}\quad\mbox{and}\quad\iw w=\frac 1m\,\I{v}{V_D}\;.
\end{equation*}
The $1/m$ factor simplifies the expressions of the linearized relative entropy and Fisher information, as we shall see below. It has no impact on the rates. Consistently with the passage to the quotient, the relative entropy and the relative Fisher information are related as follows.
\begin{proposition}\label{Prop:EntrProd} Under Assumptions {\rm (H1'')-(H2'')}, if $w$ is a solution of~\eqref{FP.Problem.Quotients}, then
\begin{equation}\label{Entr.Fisher}
\frac d{dt}\fw{w(t)}=-\,\iw{w(t)}\;.
\end{equation}
\end{proposition}
\proof As in Section~\ref{Subsec:RelMass}, consider a test
function $\phi_\lambda(x):=\chi(|x|/\lambda)$ where $\chi$ is a
smooth function on~$\RR^+$ such that $\chi\equiv 1$ on $[0,1]$,
$\chi\equiv 0$ on $[2,\infty)$, and $0\leq\chi\leq 1$ on $[1,2]$.
Then, using \eqref{FP.Problem.Quotients}, the equality
$(w^{m-1}-1)\,V_{D_*}^{m-1}=v^{m-1}- \,V_{D_*}^{m-1}$ and
integration by parts we obtain
\begin{equation*}\begin{split}
\frac{\rd}{\dt} &\frac1{1-m} \int_{\RR^d} \big[\,(w-1)-\frac{1}{m}(w^m-1)\big]\phi_\lambda\,V_{D_*}^m\dx\\
&=-\frac m{(m-1)2}\int_{\RR^d} \nabla\!\left[(w^{m-1}-1) \,V_{D_*}^{m-1}\right]\,\cdot\,\nabla\!\left[(w^{m-1}-1) \,V_{D_*}^{m-1}\,\phi_\lambda\,\right]\;w\,V_{D_*}\dx\\
&= -\frac m{(m-1)2}\int_{\RR^d}\nabla(v^{m-1}- \,V_{D_*}^{m-1})\cdot\big[\phi_\lambda\nabla(v^{m-1}- \,V_{D_*}^{m-1})+\nabla\phi_\lambda(v^{m-1}- \,V_{D_*}^{m-1})\big]\,v\dx\\
&= -\frac m{(m-1)2}\int_{\RR^d}\Big|\nabla(v^{m-1}-\,V_{D_*}^{m-1})\Big|^2\,\phi_\lambda\,v\dx + \mathcal{R}(\lambda)\;,
\end{split} \end{equation*}
where the last integral can be computed as
\begin{eqnarray*}
\mathcal{R}(\lambda) &:=& - \frac m{2\,(m-1)2}\int_{\RR^d}\nabla\!\left[(v^{m-1}-\,V_{D_*}^{m-1})2\right]\,v\,\nabla\phi_\lambda\dx\\
&&= \frac m{2\,(m-1)2}\int_{\RR^d}\left|v^{m-1}-\,V_{D_*}^{m-1}\right|^2 \big(\nabla v\cdot\nabla\phi_\lambda+v\cdot\Delta\phi_\lambda\big)\dx\;.
\end{eqnarray*}
In the region $\Omega_\lambda=\{x\in\RR^d : \lambda\le |x|\le 2\,\lambda\}$, $\lambda>0$, we get
\begin{eqnarray*}
|\,\mathcal{R}(\lambda)\,| &\le & k_1\,\int_{\Omega_\lambda}\left|v-\,V_{D_*}\right|^2\,V_{D_0}^{2(m-2)}\, \left(\,|\nabla v|\,|\nabla\phi_\lambda|+v\,|\Delta\phi_\lambda|\,\right)\dx\\
&\le & k_1\,\sup_{\Omega_\lambda}\left[V_{D_0}^{2(m-2)}\,\left|v-\,V_{D_*}\right|\,\left(\,|\nabla v|\,|\nabla\phi_\lambda|+v\,|\Delta\phi_\lambda|\,\right)\right] \,\int_{\Omega_\lambda}\left|v-\,V_{D_*}\right|\dx\;,
\end{eqnarray*}
where the positive constant $k_1$ depends on $m$, $d$, $D_0$ and $D_1$.

As in the proof of Theorem~\ref{lem:holdereg}, consider a solution $v$ of~\eqref{FP.Problem} and define
\begin{equation*}
v_\lambda(x,t)=\lambda^{2/(1-m)}\,v(\lambda\,x,t)\;.
\end{equation*}
In what follows, $c_i$ will denote positive constants which may depend on $m$, $d$, $D_0$, $D_1$ and on the maximum of $\nabla\phi_1$, but not on $\lambda$.

For any $\lambda>0$, $v_\lambda$ is a solution of~\eqref{FP.Problem} but the region $\Omega_\lambda=\{x\in\RR^d : \lambda\le |x|\le 2\lambda\}$ gets mapped into the region $\Omega_1=\{x\in\RR^d : 1\le |x|\le 2\}$. We already know that $\nabla v_\lambda$ is uniformly bounded on $\Omega_1$, by Theorem~\ref{lem:holdereg}: $\sup_{x\in\Omega_1}\big|\nabla v_\lambda(t,x)\big|\le c_0$. In terms of $v$, this gives the estimate
\begin{equation*}
\lambda^{\frac{2}{1-m}}\sup_{y\in\Omega_\lambda}\big|\nabla_y v(t,y)\big|=\lambda^{\frac{2}{1-m}}\sup_{x\in\Omega_1}\big|\lambda^{-1}\nabla_x v(t,\lambda\,x)\big|=\lambda^{-1} \sup_{x\in\Omega_1}\big|\nabla_x v_\lambda(t,x)\big|\le \lambda^{-1}\, c_0
\end{equation*}
and proves that
\begin{equation*}\label{est.grad.v}
\sup_{x\in \Omega_\lambda}\big|\nabla v(t,x)\big|\le c_0\,\lambda^{-\frac{2}{1-m}-1}\,.
\end{equation*}
By our choice of $\phi_\lambda$, we see that
\begin{equation*}\label{est.grad.lapl.varphi}
\sup_{\Omega_\lambda}|\nabla\phi_\lambda|\le \frac{c_1}{\lambda}\quad\mbox{and}\quad\sup_{\Omega_\lambda}|\Delta\phi_\lambda|\le \frac{c_2}{\lambda2}\;.
\end{equation*}
Putting together these two estimates, we get
\begin{equation*}\label{est.sum}
|\nabla v|\,|\nabla\phi_\lambda|+v\,|\Delta\phi_\lambda|\le c_3\,\lambda^{-\frac{2}{1-m}-2}\,.
\end{equation*}

Next, we observe that
\[
-\frac{\partial\,V_D}{\partial D}=\frac 1{1-m}\left[D+\frac{1-m}{2\,m}\,|x|^2 \right]^{-\frac{2-m}{1-m}}= \frac 1{1-m}\,V_D^{2-m}\;.
\]
Hence, for some constant $c_4$ depending on $m$, $d$, $D_0$ and $D_1$,
\begin{equation*}\label{est.diff.baren}
\left|V_{D_1}-\,V_{D_0}\right|\le c_4\,V_{D_0}^{2-m}
\end{equation*}
and
\[
\sup_{\Omega_\lambda}\left[V_{D_0}^{2(m-2)}\,\left|v-\,V_{D_*}\right|\,\left(\,|\nabla v|\,|\nabla\phi_\lambda|+v\,|\Delta\phi_\lambda|\,\right)\right]
\le\, c_4\,\sup_{\Omega_\lambda}\left[V_{D_0}^{m-2}\left(\,|\nabla v|\,|\nabla\phi_\lambda|+v\,|\Delta\phi_\lambda|\,\right)\right]\,.
\]
Taking into account the fact that, for any $\lambda>0$, $V_D\le\, c_5\, \lambda^{-2/(1-m)}$ on $\Omega_\lambda$, we obtain
\begin{equation*}
\sup_{\Omega_\lambda} \left[V_{D_0}^{2(m-2)}\,\left|v-\,V_{D_*}\right|\,\left(\,|\nabla v|\,|\nabla\phi_\lambda|+v\,|\Delta\phi_\lambda|\,\right)\right]\le\, c_6 \,\lambda^{2\frac{2-m}{1-m}}\,\lambda^{-\frac{2}{1-m}-2}=c_6\;,
\end{equation*}
for some positive constant $c_6$ which is independent of $\lambda$.
By assumptions (H1')-(H2') and the $\LL1$-contraction principle, the difference $v-V_{D_*}$ is in $\LL1$, and so, $\lim_{\lambda\to\infty}\int_{\Omega_\lambda}\left|v-\,V_{D_*}\right|\dx=0$. This proves that $\lim_{\lambda\to\infty}\mathcal{R}(\lambda)=0$ and we conclude by passing to the limit as $\lambda\to\infty$. \qed

\section{Convergence without rate and in relative error}\label{Sec:ConvergenceWithoutRate}

This section is mostly devoted to the proof of Theorem~\ref{Thm:A1}.

\subsection{Relative entropy}\label{sec:ferferlin}

Under Assumptions (H1'')-(H2''), the relative entropy is well defined.
\begin{lemma}[An equivalence result]\label{Lem.Bounds.RE} Let $m\in (0,1)$. If $w$ satisfies {\rm (H1'')-(H2'')}, then
\begin{equation*}
\frac 12\,W_1^{m-2}\,\int_{\RR^d}|w-1|^2\,V_{D_*}^{m}\dx \le \fw w\le \frac 12\,W_0^{m-2}\,\int_{\RR^d}|w-1|^2\,V_{D_*}^{m}\dx\,.
\end{equation*} \end{lemma}
\proof For $a>0$, let $\phi_{a}(w):=\frac{1}{1-m}\left[(w-1)-(w^m-1)/m\right] -a\left(w-1\right)2$. We compute $\phi_{a}'(w)=\frac{1}{1-m}\left[1-w^{m-1}\right]-\,2\,a\left(w-1\right)$ and $\phi_{a}''(w)=\,w^{m-2}-\,2\,a$, and note that $\phi_{a}(1)=\phi_{a}'(1)=0$. With $a=W_1^{m-2}/2$, $\phi_{a}''$ is positive on $(W_0,W_1)$, which proves the lower bound after multiplying by $V_D^m$ and integrating over $\RR^d$. With $a=\,W_0^{m-2}/2$, $\phi_{a}''$ is negative on $(W_0,W_1)$ which proves the upper bound. \qed

\begin{lemma}[Boundedness of the free energy]\label{Lem.RE.1} Let $m\in (0,1)$. If $w_0$ satisfies {\rm (H1'')-(H2'')}, then the free energy $\fw {w(t)}$ is finite for any $t\ge 0$. \end{lemma}
\proof By virtue of Proposition \ref{Prop:EntrProd}, we have to prove the result only for $w=w_0$. Notice that for any $D_0$, $D_1>0$ there exist a positive constant $C$ such that $|V_{D_0}-V_{D_1}| \le C\,|x|^{-2(\frac{2-m}{1-m})}$ as $|x|\to\infty$. Indeed,
\[
|V_{D_0}-V_{D_1}|=\left( \frac{1-m}{2\,m}\right)^{-\frac{1}{1-m}} \frac{2\,m\,|D_0-D_1|}{(1-m)2}\,|x|^{-2\left(\frac{2-m}{1-m}\right)}\,(1+o(1))\quad\mbox{as}\quad |x|\to\infty\;.
\]
By Lemma~\ref{Lem.Bounds.RE}, for some positive constant $c$ depending on $D_0$ and $D_1$, we have
\[
\frac 2{W_0^{m-2}}\,\fw w\le\int_{\RR^d}\left|\dfrac{v}{V_{D_*}}-1\right|^2V_{D_*}^m\dd x \le\int_{\RR^d}\left|v-V_{D_*}\right|^2\,V_{D_*}^{m-2}\dd x \le c\int_{\RR^d}\left|V_{D_0}-V_{D_1}\right|^2\,V_{D_*}^{m-2}\dd x\;.
\]
If $m\in (m_*,1)$, then $\left|V_{D_0}-V_{D_1}\right|^2\,V_{D_*}^{m-2}=O\big(|x|^{-\,2\, \frac{2-m}{1-m}}\big)$ is integrable as $|x|\to\infty$. Otherwise, if $m\in (0,m_*]$, $\fw w$ is also integrable as $|x|\to\infty$ because
\[
\frac 2{W_0^{m-2}}\,\fw w\le \frac 2{W_0^{m-2}}\,\fw{w_0}\le \int_{\RR^d}\left|f\right|^2\,V_{D_*}^{m-2}\dd x\le\int_{\RR^d}|f|\left|V_{D_0}-V_{D_1} \right|\,V_{D_*}^{m-2}\dd x\;,
\]
$f$ is integrable and $\left|V_{D_0}-V_{D_1}
\right|\,V_{D_*}^{m-2}$ is bounded (we ask the reader to check
this fact). \qed

\subsection{Pointwise convergence in relative error} 
\begin{lemma}\label{prop:cvrgcewithoutrate}
Let $m\in (0,1)$. If $w$ is a solution
of~\eqref{FP.Problem.Quotients} with initial data $w_0$ satisfying
{\rm (H1'')-(H2'')}, then $\lim_{t\to\infty}w(t,x)=1$ for any
$x\in\RR^d$. \end{lemma}
\proof Let $w_\tau(t,x)= w(t+\tau,x)$. By the uniform $C^k$ regularity, see Theorem~\ref{lem:holdereg}, the functions $w_\tau$ are uniformly $C1$ continuous. Hence, by the Ascoli-Arzel\`a theorem, there exists a sequence $\tau_n\to\infty$ such that $w_{\tau_n}$ converges to a function~$w_\infty$, locally uniformly in $(t,x)$. We know by the Comparison Principle, see Lemma~\ref{Lem:MP}, that $w_\infty>W_0>0$. By interior regularity of the solutions, the derivatives also converge everywhere.

By Lemma~\ref{Lem.RE.1}, $\fw w$ is finite. Since
\[
\fw {w(\tau_n)}-\fw {w(\tau_n+1)} = \int_{\tau_n}^{\tau_n+1}\iw{w(s)}\ds = \int_{0}^{1}\iw{w(s+\tau_n)}\ds\;,
\]
as a function of $t$, $\iw{w_{\tau_n}}$ is integrable on $[0,1]$ and converges to zero as $n\to\infty$,
\begin{equation*}
\lim_{n \to \infty} \int_{0}^{1}\int_{\RR^d} \left|\nabla\left[\left(w_{\tau_n}^{m-1}(t,x)-1 \right)V_{D_*}^{m-1}(x)\right]\right|^2\,w_{\tau_n}(t,x)\,V_{D_*}(x)\dx\dt=0\;.
\end{equation*}
By Fatou's lemma, $w_\infty=\lim_{n \to \infty}w_{\tau_n}$ satisfies $\nabla\left[\left(w_\infty^{m-1}-1 \right)V_{D_*}^{m-1}\right]=0$ a.e. in $(0,1)\times\RR^d$. As a consequence of the conservation of relative mass, see Proposition~\ref{prop:relconsmass}, $w_\infty=1$. Thus, we have proved the convergence a.e., and by equi-continuity, the pointwise convergence. Since the limit is unique, the whole family $\{w_\tau\}_\tau$ converges everywhere as $\tau\to\infty$. \qed

\subsection{Proof of Theorem~\ref{Thm:A1}}

\subsubsection*{Proof of Theorem~\ref{Thm:A1}, (i) and (ii)}

By Lemma~\ref{prop:cvrgcewithoutrate}, $\lim_{t \to \infty}\left|v(t,x)-V_{D_*}(x)\right|=0$ for any $x\in\RR^d$. Moreover, we observe that
\[
\left|v(t)-V_{D_*}\right|\leq\max\big\{\left|V_{D_0}-V_{D_*}\right|,\; \left|V_{D_1}-V_{D_*}\right|\big\}=O\left(|x|^{-2(2-m)/(1-m)}\right)
\]
as $|x|\to\infty$. By Lebesgue's dominated convergence theorem, $v(t)$ converges to $V_{D_*}$ in $\LL^p(\RR^d)$, for any $p\in (\q(m),\infty)$, where $\q(m):= \tfrac{d\,(1-m)}{2\,(2-m)}$ is the infimum of all positive $p$ such that the difference between two different Barenblatt profiles belongs to $\LL^p(\RR^d)$.

The uniform convergence is based on the following interpolation
lemma, due to Nirenberg, \cf~\cite[p. 126]{MR0109940}. Let
$\lambda$, $\mu$ and $\nu$ be such that $-\infty<\lambda \le \mu
\le \nu <\infty$. Then there exists a positive constant $\mathcal
C_{\lambda,\,\mu,\,\nu}$ such that
\begin{equation}\label{eq:interpolation}
\|f\|_{1/\mu}^{\nu-\lambda} \le \mathcal
C_{\lambda,\,\mu,\,\nu}\,\|f\|_{1/\lambda}^{\nu-\mu} \;
\|f\|_{1/\nu}^{\mu-\lambda}\quad\forall\;f\in\mathcal C(\RR^d)\;,
\end{equation}
where $\|\cdot\|_{1/\sigma}$ stands for the following quantities:

\noindent(i) If $\sigma>0$, then
$\|f\|_{1/\sigma}=\left(\int_{\RR^d}|f|^{1/\sigma}\dx\right)^\sigma$.

\noindent(ii) If $\sigma<0$, let $k$ be the integer part of
$(-\sigma\,d)$ and $\alpha=|\sigma|\,d-k$ be the fractional
(positive) part of $\sigma$. Using the standard multi-index
notation, where $|\eta|=\,\eta_1+\,\ldots\,+\eta_d$ is the length
of the multi-index
$\eta=(\eta_1,\,\ldots\,\eta_d)\in\mathbb{Z}^d$, we define
\begin{equation*}\label{def.C^k}
\|f\|_{1/\sigma}=\left\{\begin{array}{lll} \displaystyle\max_{|\eta|=k}\; \big|\partial^\eta f\big|_\alpha=\displaystyle\max_{|\eta|=k}\; \sup_{x,\,y\in\RR^d}\;\dfrac{\big|\partial^\eta\,f(x)-\partial^\eta\,f(y)\big|}{|x-y|^\alpha}=:\|f\|_{C^\alpha(\RR^d)}& \mbox{if~}\alpha>0\;,\\[5mm]
\displaystyle\max_{|\eta|=k}\;\displaystyle\sup_{z\in\RR^d}\big|\partial^\eta f(z)\,\big|\,:=\|f\|_{C^k(\RR^d)}& \mbox{if~}\alpha=0\;.
\end{array}\right.
\end{equation*}
As a special case, we observe that $\|f\|_{-d/j}=\|f\|_{C^{j}(\RR^d)}$.

\noindent(iii) By convention, we note $\|f\|_{1/0}=\sup_{z\in\RR^d}|f(z)|=\|f\|_{C0(\RR^d)}=\|f\|_{\infty}$.

Let $j\in\mathbb{N}$ and $\lambda=-(j+1)/d\le\mu=-j/d\le\nu=1/2$ so that $k=j+1$ and $\alpha=0$. Inequality~\eqref{eq:interpolation} becomes
\begin{equation}\label{interp.C_j}
\|f\|_{C^{j}(\RR^d)}\;\le\;\mathcal C_{-(j+1)/d,\,-j/d,\,1/2}^\frac{2d}{d+2(j+1)}\;\|f\|_{C^{j+1}(\RR^d)}^{\frac{d+2j}{d+2(j+1)}}\;\|f\|_2^{\frac{2}{d+2(j+1)}}
\end{equation}
for any $j\in\mathbb{N}$. By applying this interpolation inequality $f=v(t)-V_{D_*}$ with $j=0$, we obtain
\begin{equation}\label{interp.C_j2}
\|v(t)-V_{D_*}\|_{\infty}\;\le\;\mathcal C_{-1/d,\,0,\,1/2}^\frac{2d}{d+2}\;\|v(t)-V_{D_*}\|_{C^{1}(\RR^d)}^{\frac{d}{d+2}}\;\|v(t)-V_{D_*}\|_2^{\frac{2}{d+2}}\;.
\end{equation}
By Theorem~\ref{lem:holdereg}, the $C1$ norm is uniformly bounded. If $q(m)<2$, that is, if $d\le 8$, or $d\ge 9$ and $m>(d-8)/(d-4)$, we already know that $\lim_{t \to \infty}\|v(t)-V_{D_*}\|_2=0$. Otherwise, we can interpolate $\|v(t)-V_{D_*}\|_2$ between $\|v(t)-V_{D_*}\|_1\leq\|v_0-V_{D_*}\|_1$ (see Lemma~\ref{Lem:contraction}) and $\|v(t)-V_{D_*}\|_q$ for some $q>q(m)$. This proves that $\lim_{t \to \infty}\|v(t)-V_{D_*}\|_\infty=0$.\qed

\subsubsection*{Proof of Theorem~{\rm \ref{Thm:A1}}, (iii)}
\begin{corollary}[Uniform convergence of the relative error]\label{cor3.4} Let $m\in (0,1)$. If $w$ is a solution of~\eqref{FP.Problem.Quotients} with initial data $w_0$ satisfying {\rm (H1'')-(H2'')}, then
\[
\lim_{t\to\infty}\|w(t)-1\|_{\infty}=0\;.
\]\end{corollary}
\proof Because of the convergence of $v(t)$ to $V_{D_*}$ in $\LL^\infty(\RR^d)$, we know that $w(t)$ converges uniformly to $1$ on any compact set of $\RR^d$. By Assumption (H1'), $v(t)$ is sandwiched between two Barenblatt profiles that have the same asymptotic behavior when $|x|$ is large. In terms of $w$, this means that $|w(t,x)-1|$ is small for $|x|$ large, uniformly in $t$. Global uniform convergence follows.\qed

The fact that $w(t)$ converges uniformly to $1$ as $t\to\infty$ allows us to improve the lower and upper bounds $W_0$ and $W_1$ for the function $w(t)$, at the price of waiting some time. For any $\varepsilon>0$ there exists a time $t_0=t_0(\varepsilon)\ge 0$ such that
\begin{equation*}
1-\varepsilon\le w(t,x)\le 1+\varepsilon\quad\forall\;(t,x)\in (t_0,\infty)\times\RR^d\;.
\end{equation*}

\begin{corollary}[$\LL^p$ Convergence of the relative error]\label{cor3.4bis} Let $m\in (0,1)$. If $w$ is a solution of~\eqref{FP.Problem.Quotients} with initial data $w_0$ satisfying {\rm (H1'')-(H2'')}, then $w(t)$ converges to $1$ in $\LL^p(\RR^d)$ for any $p\in(d/2,\infty]$. \end{corollary}
\proof By Assumptions (H1'')-(H2''), there exists a positive constant $c_0$ such that $w_0-1$ is bounded and for $|x|$ large,
\[
|w_0-1|=\left|\frac{v_0-V_{D_*}}{V_{D_*}}\right|\le\frac{V_{D_1}-V_{D_0}}{V_{D_*}}\le \frac{c_0}{1+|x|^2}\;.
\]
By Lemma~\ref{Lem:MP}, the same estimate holds for $w(t)$. Hence $w(t)-1\in\LL^q$ for any $q>d/2$. Let $\delta=(p-d/2)/2>0$. By H\"older's inequality,
\[
\lim_{t \to \infty}\int_{\RR^d}\left|w(t)-1\right|^{p}\dx\le\lim_{t \to \infty}\left\|w(t)-1\right\|_{\infty}^\delta\int_{\RR^d}\left(\frac{c_0}{1+|x|^2}\right)^{\delta+d/2}\!\!\dx=0\;.
\]
\qed

Theorem~\ref{Thm:A1}, (iii), results from Corollaries~\ref{cor3.4} and \ref{cor3.4bis}.

\subsection{Uniform convergence and $C^\alpha$ regularity}
\begin{lemma}\label{Lemma:Calpha} Let $m\in (0,1)$. Consider a solution $v$ of~\eqref{FP.Problem} with initial data $v_0$ satisfying {\rm (H1')-(H2')}. There exists $t_0\geq 0$, $\alpha \in (0,1)$ and a positive constant $\mathcal{H}$ such that $h(t):=v(t)-V_{D_*}$ is in $C^\alpha$ and
\begin{equation}\label{eq:mathcalH}
\|h(t)\|_{C^\alpha(\RR^d)}\le\,\mathcal{H}\,\|h(t)\|_\infty\quad\forall\; t\geq t_0\;.
\end{equation} \end{lemma}
\proof Since both $v$ and $V_{D_*}$ are solutions to equation \eqref{FP.Problem}, $h$ solves
\begin{equation*}\label{Eq.Difference}
h_t=\nabla\cdot\left[m(V_{D_*}+h)^{m-1}\nabla\,h +m\left((V_{D_*}+h)^{m-1}-V_{D_*}^{m-1}-V_{D_*}^{m-2}h\right)\nabla\,V_{D_*} \right]\,.
\end{equation*}
Let $\lambda>0$. By Theorem~\ref{Thm:A1}, we know that for some $t_0\geq 0$, for any $t\geq t_0$, $\|h(t)\|_\infty$ can be taken uniformly small and $v$ uniformly positive on $B_{2\lambda}$. We apply the standard quasilinear parabolic theory, see {\sl e.g.\/}~\cite[Theorem~1.1, p.~418]{LSU}, with structure functions $a_i(x,t,h,\xi)= A\,\xi+ B\,h$ and $a=0$, where $A(x,t):=m\,v^{m-1}$, $B(x,t):=m\,\big[\big(\frac{v^{m-1}-V_{D_*}^{m-1}} {v-V_{D_*}}-V_{D_*}^{m-2} \big)\,\nabla\,V_{D_*}\big]$. Hence there exists a H\"older exponent $\alpha\in (0,1)$ and a constant $\mathcal{H}$ depending on the uniform bounds for the coefficients, and on $\lambda$, such that~\eqref{eq:mathcalH} is verified in $B_\lambda\times (t_0+1,\infty)$. To extend the estimate uniformly to the whole space, $x\in \RR^d$, we use the same scaling argument as in the proof of Theorem~\ref{lem:holdereg}. We leave the details to the reader.\qed

\section{Linearization} \label{Sec:SpectralGap}

In order to better understand the asymptotic behavior of the solutions of \eqref{FP.Problem.Quotients}, we linearize the equation around the equilibrium, introducing a convenient weight. Let $g$ be such that
\begin{equation}\label{g.Linearization}
w(t,x)=1+\varepsilon\,\frac{g(t,x)}{V_{D_*}^{m-1}(x)}\quad\forall\; t>0\;,\quad\forall\;x\in\RR^d\;,
\end{equation}
for some $\varepsilon>0$, small. Substituting this expression in Equation~\eqref{FP.Problem.Quotients} and letting $\varepsilon\to 0$, we formally obtain a linear equation for $g$,
\begin{equation} \label{Linearised.FP}
g_t=A_m\,g\quad \mbox{where}\quad A_m\,g:=m\,V_{D_*}^{m-2}(x)\,\nabla\cdot\left[V_{D_*}\,\nabla g\right]\;.
\end{equation}
The linear operator
$A_m:\LL2(\RR^d,V_{D_*}^{2-m}\dx)\to\LL2(\RR^d,V_{D_*}^{2-m}\dx)$
is the positive self-adjoint operator associated to the closure of
the quadratic form defined for $\phi\in C_c^{\infty}(\RR^d)$ by
\begin{equation}\label{Dirichlet.Form} \IL\phi:=m\int_{\RR^d}|\nabla g|^2\,V_{D_*}\dx\;. \end{equation} See~\cite[Theorem 2.6]{MR990239} for more details.

\medskip

With the same heuristics, we linearize the relative entropy
$\mathcal F$ and the relative Fisher information~$\mathcal J$,
which provides the functionals $\textsf F$ and $\textsf I$, where
$\textsf I$ is given by \eqref{Dirichlet.Form} and $\textsf F$ is
defined by
\begin{equation}
\EL g:=\frac{1}{2} \int_{\RR^d}|g|^2\,V_{D_*}^{2-m}\dx\;.
\end{equation}
Note that $\EL g$ is the $\LL2(\RR^d,V_{D_*}^{2-m}\dx)$-norm up to a factor $1/2$. If $g$ is a solution of~\eqref{Linearised.FP}, then
\begin{equation}
\frac{\rd}{\dt}\EL{g(t)}=- \,\IL{g(t)}\;.
\end{equation}
Still at a formal level, the conservation of relative mass amounts to require
\begin{equation*}
\int_{\RR^d}(v_0-V_{D_*})\dx=\int_{\RR^d}(w-1)\,V_{D_*}\dx= \varepsilon \int_{\RR^d}g\,V_{D_*}^{2-m}\dx
\end{equation*}
in the limit $\varepsilon \to 0$. Hence, it makes sense to require that $\int_{\RR^d}g\,V_{D_*}^{2-m}\dx=0$ and use the spectral gap estimate, see~\cite{BBDGV-CRAS} and Theorem~\ref{Teor.Spect.}. With $\Cmd=m/\gopt$, we obtain
\begin{equation}\label{Spectral.Gap}
  2 \,\EL g\le \frac{\Cmd}m\,\IL g\;,
\end{equation}
which gives, for the solution of~\eqref{g.Linearization}, an
exponential decay of the relative entropy,
\begin{equation*}\label{Exp.Decay.Linear}
\EL{g(t)}\le {\rm e}^{-\,2\,\gopt\,t}\,\EL{g(0)}\quad\forall\; t\ge 0\;.
\end{equation*}

\medskip In Sections \ref{sec:cvrgcewithoutrate} and \ref{Sec:WithRate}, we will compare the relative entropy estimates for the solutions of~\eqref{FP.Problem.Quotients} with the ones of the linearized problem. This is the main ingredient of the proof of Theorem~\ref{Thm:A2}.

\medskip The connection with the Fokker-Planck equation is easy to understand at the level of the linearized problem. In the limit $m\to 1$, we observe that
\begin{equation*}
\lim_{m\to 1^-}D_*^{1/(1-m)}\,V_{D_*}=(2\,\pi\,D_*)^{d/2}\,\mu \quad\mbox{with}\quad\mu(x)=\frac{{\rm e}^{-\frac{|x|^2}{2\,D_*}}}{(2\,\pi\,D_*)^{d/2}}\;.
\end{equation*}
Equation~\eqref{Linearised.FP} formally converges to the Ornstein-Uhlenbeck equation,
\begin{equation*}\label{Orn.Uhl}
g_t=\mu^{-1}\,\nabla\cdot\big(\,\mu\,\nabla g\,\big)\;.
\end{equation*}
The spectral gap inequality~\eqref{Spectral.Gap} corresponds in such a limit to the well-known Poincar\'e inequality with gaussian weight,
\begin{equation*}\label{Poincare.Gauss}
\int_{\RR^d}|\phi|^2\,\rd\mu\le\int_{\RR^d}\left|\nabla \phi\right|^2\,\rd\mu\quad\forall\;\phi\in C^\infty(\RR^d)\;\mbox{such that}\int_{\RR^d}\phi\;\rd\mu=0\;,
\end{equation*}
where $\rd\mu:=\mu\dx$. Note that in the Gaussian case, a logarithmic Sobolev inequality holds, see~\cite{Gr-LSI},
\begin{equation*}\label{Gross.LSI}
\int_{\RR^d}|\phi|^2\,\log\left(\frac{|\phi|^2}{\int_{\RR^d}|\phi|^2\,\rd\mu}\right) \;\rd\mu\;\le\;2\int_{\RR^d}\left|\nabla \phi\right|^2\,\rd\mu\;,
\end{equation*}
which is stronger than the Gaussian Poincar\'e inequality. This is not the case with the measure $V_{D_*}\!\dx$. Although the spectral gap inequality~\eqref{Spectral.Gap} holds true, there is no corresponding logarithmic Sobolev inequality.

\section{More on the relative Fisher information}\label{sec:cvrgcewithoutrate}

In this section, we relate the relative Fisher and linearized Fisher informations. This and Lemma~\ref{Lem.Bounds.RE} provide us with an estimate of the relative entropy in terms of the relative Fisher information, or {\sl entropy - entropy production inequality,\/} for the nonlinear problem.

\subsection{Fisher information and linearized Fisher information}\label{sec:fishfishlin}

\begin{lemma}[Upper bound on the Fisher information]\label{Lem.Bounds.fish} Let $m\in (0,1)$. There exists two positive constants $\beta_1$ and $\beta_2$ (depending on $W_0$, $W_1$ and $m$) such that, for any $w$ satisfying {\rm (H1'')-(H2'')},
\begin{equation*}
\IL g\le\beta_1\,\iw w+\beta_2\,\EL g\quad\mbox{with}\quad g:=(w-1)\,V_{D_*}^{m-1}\,.
\end{equation*}
Moreover, if $\eta:=\max\{1-W_0,W_1-1\}$, then $\lim_{\eta\to 0_+}\big(|1-\beta_1|+\beta_2\big)=0$. \end{lemma}
The constant $\beta_1$ and $\beta_2$ are explicitly given in \eqref{eq:c0c1c2} in terms of $m$, $D_*$, $W_0$ and $W_1$.

\proof Define $h_k(w):=(w^{k-1}-1)/(k-1)$. Let $\alpha_0:=W_0^{2(2-m)}$, $\alpha_1:=W_1^{2(2-m)}$. Since $|h_2/h_m|$ is non-decreasing,
\begin{equation}\label{eq:c0c1}
\alpha_0 \le\left|\dfrac{h_2'(W_0)}{h_m'(W_0)}\right|^2\le\left|\dfrac{h_2(w)} {h_m(w)}\right|^2\le\left|\dfrac{h_2(W_1)}{h_m(W_1)}\right|^2\le\left| \dfrac{h_2'(W_1)}{h_m'(W_1)}\right|^2=\alpha_1\;.
\end{equation}
Note that $\alpha_0=\alpha_0(W_0)<1<\alpha_1=\alpha_1(W_1)$ and both converge to $1$ as $W_0, W_1\to 1$.

\medskip Using the fact that $V_{D_*}^{m-1}=D_*+\frac{1-m}{2\,m}\,|x|^2$ and an integration by part, we~get
\begin{equation*} \begin{split}
&\int_{\RR^d}\left|\nabla\left[h_k(w)V_{D_*}^{m-1}\right]\right|^2\,V_{D_*}\dx\\
&= \int_{\RR^d}|h_k'(w)|^2\,|\nabla w|^2\,V_{D_*}^{2m-1} \dx + \frac{1-m}{m2} \!\int_{\RR^d}|x|^2\,|h_k(w)|^2\,V_{D_*} \dx -d\,\frac{1-m}{m}\!\int_{\RR^d}|h_k(w)|^2\,V_{D_*}^m \dx\;.
\end{split} \end{equation*}
Let $g:=(w-1)\,V_{D_*}^{m-1}$. Applied with $k=2$ and $k=m$, the above identity gives
\begin{multline*}
\IL g= m\,\int_{\RR^d}\left|\nabla\left[h_2(w)V_{D_*}^{m-1}\right]\right|^2\,V_{D_*}\dx\\
\hspace*{-4cm}\le m\,\alpha_1\,\int_{\RR^d}|h_m'(w)|^2\,|\nabla w|^2\,V_{D_*}^{2m-1} \dx \\
+ \alpha_1\,\frac{1-m}{m} \int_{\RR^d}|x|^2\,|h_m(w)|^2\,V_{D_*} \dx -d\,(1-m)\int_{\RR^d}|h_2(w)|^2\,V_{D_*}^m \dx
\end{multline*}
and
\begin{multline*}
\int_{\RR^d}|h_m'(w)|^2\,|\nabla w|^2\,V_{D_*}^{2m-1} \dx\\
\hspace*{-4cm}=\int_{\RR^d}\left|\nabla\left[h_m(w)V_{D_*}^{m-1}\right]\right|^2\,V_{D_*}\dx\\
- \frac{1-m}{m2} \int_{\RR^d}|x|^2\,|h_m(w)|^2\,V_{D_*} \dx +d\,\frac{1-m}{m}\int_{\RR^d}|h_m(w)|^2\,V_{D_*}^m \dx\;.
\end{multline*}
Collecting these estimates, we obtain
\begin{equation*}
\IL g\le m\,\alpha_1\,\int_{\RR^d}\left|\nabla\left[h_m(w)V_{D_*}^{m-1}\right]\right|^2\,V_{D_*}\dx + d\,(1-m) \int_{\RR^d}\left( \alpha_1\,|h_m(w)|^2 - |h_2(w)|^2\right) V_{D_*}^m \dx\;.
\end{equation*}
Note that
\begin{equation*}
m \int_{\RR^d}\left|\nabla\left[h_m(w)V_{D_*}^{m-1}\right]\right|^2\,V_{D_*}\dx \le \,W_0^{-1}\,\iw w\;.
\end{equation*}
Using $\EL g= \frac 12\int_{\RR^d}|g|^2\,V_{D_*}^{2-m} \dx=\frac 12\int_{\RR^d}|h_2(w)|^2\,V_{D_*}^m \dx$ with $g:=(w-1)\,V_{D_*}^{m-1}$, we obtain
\begin{equation*}
\IL g\le\beta_1\,\iw w + \beta_2\,\EL w\;,
\end{equation*}
with
\begin{equation}\label{eq:c0c1c2}
\beta_1:=\frac{\alpha_1}{W_0}=\frac{W_1^{2(2-m)}}{W_0}\quad \mbox{and}\quad \beta_2:=2\,d\,(1-m)\left(\frac{\alpha_1}{\alpha_0} -1\right) \;.
\end{equation}
Note that $\alpha_0=\alpha_0(W_0)<1<\alpha_1=\alpha_1(W_1)$ and
both tend to 1 as $W_0, W_1\to 1$. \qed

\subsection{Entropy - entropy production inequality}
\begin{theorem}[Entropy - entropy production inequality]\label{Cor:E-I} Let $m\in (0,1)$, $m\neq m_*$. For any function $w=w(x)$ satisfying {\rm (H1'')-(H2'')}, if $1-W_0>0$ and $W_1-1>0$ are small enough, then there exists a positive constant $\gamma$ such that
\begin{equation*}\label{eq:gammapasopt}
\gamma\,\fw w\leq\iw w\;.
\end{equation*} \end{theorem}
As we shall see in the proof, the constant $\gamma$ can be estimated as follows:
\begin{equation}\label{Eqn:ExprGamma}
\gamma\geq 2\,\frac{m-\Cmd\,d\,(1-m)\left[ \left( \frac{W_1}{W_0}\right)^{2(2-m)}-1\right]}{\Cmd\,W_0^{m-3}\,W_1^{2(2-m)}}\;.
\end{equation}
The condition that {\sl $1-W_0>0$ and $W_1-1>0$ are small enough\/} in the statement of Theorem~\ref{Cor:E-I} can be replaced by a weaker condition which amounts to ask that the right hand side of the above estimate is positive, that is, with $\gopt=m/\Cmd$,
\begin{equation*}\label{Explicit}
\frac{W_1}{W_0}<\left(1+\frac{\gopt}{d\,(1-m)}\right)^\frac 1{2\,(2-m)}\,.
\end{equation*}
\proof Let $g:=(w-1)\,V_{D_*}^{m-1}$. By definition of $D_*$, $0=\int_{\RR^d}(w-1)\,V_{D_*}\dx=\int_{\RR^d} g\,V_{D_*}^{2-m}\dx$ if $m>m_*$. By the spectral gap estimate~\eqref{Spectral.Gap} (also see Theorem~\ref{Teor.Spect.}),
\begin{equation*}
2 \,\EL g\le \frac{\Cmd}m\,\IL g\;.
\end{equation*}
By Lemma~\ref{Lem.Bounds.fish},
\[
2\,\EL g\le\frac{\Cmd}m\,\IL g\le\frac{\Cmd}m\,(\beta_1\,\iw w+\beta_2\,\EL g)\;,
\]
from which we deduce that
\[
\EL g\le \frac{\Cmd\,\beta_1}{2\,m-\Cmd\,\beta_2}\,\iw w\;.
\]
By Lemma~\ref{Lem.Bounds.RE}, the conclusion holds with
\begin{equation*}\label{def.gamma}
\gamma=\frac{2\,m-\Cmd\,\beta_2}{W_0^{m-2}\,\Cmd\,\beta_1}
\end{equation*}
under the condition $2\,m-\Cmd\,\beta_2>0$. According to the definition~\eqref{eq:c0c1c2} of $\beta_2$, this amounts to the condition
\begin{equation*}\label{Cond.alpha}
\Cmd\,\beta_2=2\,d\,(1-m)\,\Cmd\,\left(\frac{\alpha_1}{\alpha_0}-1\right)<2\,m\;,
\end{equation*}
that is $\alpha_1/\alpha_0$ close enough to $1$, which follows from the requirement that $1-W_0>0$ and $W_1-1>0$ are small enough.\qed

\section{Convergence with rates}\label{Sec:WithRate}

\subsection{Exponential decay of the relative entropy}

\begin{proposition}\label{Thm:expodecfreeenergy} Let $m\in (0,1)$, $m\ne m_*$. There exists a positive constant~$\gamma$ such that, for any solution $w$ of~\eqref{FP.Problem.Quotients} with initial data $w_0$ satisfying~{\rm (H1'')-(H2'')}, if $1-W_0>0$ and $W_1-1>0$ are small enough, then
\begin{equation*}\label{Exp.Decay.Entropy}
\fw {w(t)}\le \fw{w_0}\;\ee^{-\gamma\,t}\;.
\end{equation*} \end{proposition}
The value of $\gamma$ can be estimated from below by \eqref{Eqn:ExprGamma}.
\proof We combine formula $\frac d{dt}\fw {w(t)}=-\iw {w(t)}$ with the estimate $\iw {w(t)}\ge \gamma\,\fw {w(t)}$ obtained in Theorem~\ref{Cor:E-I}, and then integrate the resulting differential inequality.\qed

\subsection{Moments, $\LL^p$ and $C^k$ estimates}\label{sec:interpolation}

We recall that $q_*:=\dfrac{2\,d\,(1-m)}{d\,(1-m)+2\,(2-m)}$ and define
\begin{equation*}\label{gamma.q}
\gamma(q):=\frac 12\,\;\mbox{if}\; q\in(q_*,2]\;, \quad\gamma(q)=\frac{q+d}{q(d+2)}\,\;\mbox{if}\; q\in(2,\infty)\;, \quad\gamma(\infty):=\frac 1{d+2}\;.
\end{equation*}
The following lemma helps to understand better the consequences of the convergence of the free energy $\fer{v}{V_{D_*}}{}\,$, in terms of $\LL^p$, moment and also $C^k$ convergence.
\begin{lemma}\label{Conv.Weight} Let $m\in (0,1)$ and consider a function $v$ satisfying {\rm (H1')-(H2')}. Then
\begin{enumerate}
\item[{\rm (i)}] For any $\vartheta\in[0,\frac{2-m}{1-m}]$, there exists a positive constant $K_\vartheta$ such that
\begin{equation*}
\left\|\,|x|^\vartheta(v-V_{D_*})\,\right\|_{2}\leq K_\vartheta\,\left( \fer {v}{V_{D_*}}\,\right)^{1/2}\;.
\end{equation*}
\item[{\rm (ii)}] For any $q\in(q_*,2]$, there exists a positive constant $K(q)$ such that
\begin{equation*}
\|v-V_{D_*}\|_q\le K(q)\,\left( \fer{v}{V_{D_*}}\,\right)^{\gamma(q)}\;.
\end{equation*}
\end{enumerate}
Consider now a solution $v$ of \eqref{FP.Problem} such that $v_0$ satisfies {\rm (H1')-(H2')} and fix some $t_0>0$. Then
\begin{enumerate}
\item[{\rm (iii)}] For any $j\in\mathbb{N}$ and any $t_0>0$, there exists a positive constant $H_j$ such that
\begin{equation*}
\|v(t)-V_{D_*}\|_{C^j(\RR^d)}\le H_j\,\left( \fer{v(t)}{V_{D_*}}\,\right)^{\frac{1}{d+2(j+1)}}\quad\forall\;t\geq t_0\;.
\end{equation*}
\item[{\rm (iv)}] For any $q\in(2,\infty]$, there exists a positive constant $K(q)$ such that
\begin{equation*}
\|v(t)-V_{D_*}\|_q\le K(q)\,\left( \fer{v(t)}{V_{D_*}}\,\right)^{\gamma(q)}\quad\forall\;t\geq t_0\;.
\end{equation*}
\end{enumerate} \end{lemma}
\proof (i) With $\kappa_\vartheta:=2\,\sup_{r>0}r^{2\vartheta}\left(D_*+\frac{1-m}{2\,m}\,r2\right)^{-(2-m)/(1-m)}$,
\begin{equation*}
\|\,|x|^\vartheta(v-V_{D_*})\,\|_{2}2 \leq\frac 12\,\kappa_\vartheta\int_{\RR^d}|v-V_{D_*}|^2\,V_{D_*}^{m-2}\dx\le \kappa_\vartheta\,\EL{(w-1)V_{D_*}^{m-1}}\;,
\end{equation*}
and the right hand side is equivalent to $\fer {v}{V_{D_*}}$ by Lemma~\ref{Lem.Bounds.RE}.

\medskip\noindent(ii) Let $q\in(q_*,2)$. By H\"older's inequality,
\[
\int_{\RR^d}|v-V_{D_*}|^q\dx=\int_{\RR^d}V_{D_*}^{(2-m)q/2}\cdot\left(|v-V_{D_*}|^2\,V_{D_*}^{m-2}\right)^{q/2}\dx\leq c(q)\left(\int_{\RR^d}|v-V_{D_*}|^2\,V_{D_*}^{m-2}\dx\right)^{q/2}\,,
\]
where $c(q)^{2/(2-q)}:=\int_{\RR^d}V_{D_*}^{(2-m)\frac q{2-q}}\dx$ is finite for any $q>q_*$. By Lemma~\ref{Lem.Bounds.RE}, the estimate holds with $K(q):=c(q)^{1/q}\,(2\,W_1^{2-m})^{1/2}$. In the limit case $q=2$, the same method applies with $c(2)=\|V_{D_*}^{(2-m)}\|_{\infty}=D_*^{-(2-m)/(1-m)}$.

\medskip\noindent(iii) We apply the interpolation inequality~\eqref{interp.C_j} to $f=v(t)-V_{D_*}$ and bound $\|v(t)-V_{D_*}\|_{C^{j}(\RR^d)}$ in terms of $\|v(t)\|_{C^{j+1}}$, which is uniformly bounded by Theorem~\ref{lem:holdereg} and $\|v(t)-V_{D_*}\|_2^{2/(d+2(j+1))}$, for which we apply the result of Part (i) with $\vartheta=0$.

\medskip\noindent(iv) By Theorem~\ref{lem:holdereg}, $v(t)\in C1(\RR^d)$ and $v(t)-V_{D_*}$ is bounded in $C1$ uniformly for any $t\geq t_0>0$. By~\eqref{interp.C_j2},
\[
\|v(t)-V_{D_*}\|_{\infty}\;\le\;\mathcal C_{-1/d,\,0,\,1/2}^\frac{2d}{d+2}\;\|v(t)-V_{D_*}\|_{C^{1}(\RR^d)}^{\frac{d}{d+2}}\;\|v(t)-V_{D_*}\|_2^{\frac{2}{d+2}}\;.
\]
We conclude by using H\"older's inequality, $\|v(t)-V_{D_*}\|_q\le \|v(t)-V_{D_*}\|_{\infty}^{(q-2)/q}\,\|v(t)-V_{D_*}\|_{2}^{2/q}$.
\qed

\subsection{Improvement of the convergence}\label{sec.proof.CRE-exp}

\begin{theorem}\label{Thm:6.3} Let $d\ge3$, $m\in(0,1)$ with $m \ne m^*$. Consider a solution $w$ of~\eqref{FP.Problem.Quotients} with initial data satisfying {\rm (H1'')-(H2'')}. There exist a positive constant $\mathcal K$ and a time $t_0\geq 0$ such that
\begin{equation*}
\fw {w(t)}\le\,\mathcal{K}\;\ee^{-2\,\gopt\,t}\quad\forall\;t\ge t_0\;.
\end{equation*}
Moreover, for any $\lambda\in\big(0,\gopt\big)$, there exist a positive constant $\mathcal C_\infty$ and a time $t_0\geq 0$ such that
\begin{equation*}
\|w(t)-1\|_\infty\le\,\mathcal{C}_\infty\;\ee^{-2\,\frac{1-m}{2-m}\,\frac\lambda{d+2}\,t}\quad\forall\;t\ge t_0\;.
\end{equation*}\end{theorem}
Here $\gopt=m/\Cmd$ where $\Cmd$ is given in Theorem~\ref{Teor.Spect.}. Hence the rate of decay obtained by spectral methods for the linearized equation exactly gives the rate of decay for the nonlinear problem, and the price to be paid is only on the constant $\mathcal K$. As a subproduct of the proof, for some positive constants $\eta_0$ and $\gamma_\infty$ which are defined below, we obtain the following estimate
\begin{equation*}\label{EstimMathcalK}
\mathcal K\leq\fw{w(t_0)}\,\ee^\frac{\eta_0\,(2-m)}{\gamma_\infty\,(1-m)}\,\ee^{2\,\gopt\,t_0}\;.
\end{equation*}
\proof By Corollary~\ref{cor3.4}, for any $\varepsilon>0$, there exists $t_0>0$ such that $\tilde w(t)=w(t+t_0)$ satisfies Assumption (H2'') at $t=0$ with $0<1-W_0<\varepsilon$ and $0<W_1-1<\varepsilon$. With $\varepsilon$ small enough, $\tilde w$ enters in the framework of Proposition~\ref{Thm:expodecfreeenergy}, with $\gamma$ as in Theorem~\ref{Cor:E-I}. From now on, we assume that $t\geq t_0$ and simply write $w$ instead of $\tilde w$.

On the one hand, by Lemma~\ref{Conv.Weight} and Proposition~\ref{Thm:expodecfreeenergy}, we have
\begin{equation*}
\|v(t)-V_{D_*}\|_\infty \le \sigma_0\,\ee^{-\gamma_\infty\,(t-t_0)} \quad \mbox{with}\quad \sigma_0:=K(\infty)\,\fw {w_0}^{\frac 1{d+2}}\quad \mbox{and}\quad \gamma_\infty=\frac{\gamma}{d+2}\;,
\end{equation*}
which, in terms of $w$, gives the estimate
\begin{equation*}\label{Up.Error.Estimates}
|w(t,x)-1|\le \sigma_0\,\ee^{-\gamma_\infty\,(t-t_0)}\left[D_*+\frac{1-m}{2\,m}\,|x|^2\right]^{\frac{1}{1-m}}\quad\forall\;t\ge t_0\;,\quad\forall\;x\in\RR^d\;.
\end{equation*}

On the other hand, let $h_\alpha(s):=(1+s)^\alpha$, $\alpha>1$. For any $s\in[0,s_0]$,
\[
\frac{h_\alpha(s)-1}s\leq\alpha+\frac{s_0}2\,\max_{s\in[0,s_0]}h_\alpha''(s)=\left\{\begin{array}{ll}
\alpha+\frac{s_0}2\,\alpha\,(\alpha-1)\,(1+s_0)^{\alpha-2}\quad&\mbox{if}\quad\alpha\geq 2\;,\vspace*{12pt}\\
\alpha+\frac{s_0}2\,\alpha\,(\alpha-1)\quad&\mbox{if}\quad\alpha\leq 2\;.
\end{array}\right.
\]
Apply then this inequality with $\alpha=1/(1-m)$, $s=(D_*-D_1)/\big(D_1+\frac{1-m}{2\,m}\,|x|^2\big)\leq s_0=\frac{D_*}{D_1}-1$ to get the existence of a positive constant $\mathcal M_1=\mathcal M_1(m,D_*,D_1)$ such that
\[
w(t,x)-1\le\frac{V_{D_1}}{V_{D_*}}-1=\left(1+\frac{D_*-D_1}{D_1+\frac{1-m}{2\,m}\,|x|^2}\right)^\frac 1{1-m}-1\le\frac{\mathcal M_1}{D_1+\frac{1-m}{2\,m}\,|x|^2}\quad\forall\; x\in\RR^d\;.
\]
Similarly, for any $s\in[-s_0,0]$,
\[
\frac{h_\alpha(s)-1}s\geq\alpha-\frac{s_0}2\,\max_{s\in[-s_0,0]}h_\alpha''(s)=\left\{\begin{array}{ll}
\alpha-\frac{s_0}2\,\alpha\,(\alpha-1)\quad&\mbox{if}\quad\alpha\geq 2\;,\vspace*{12pt}\\
\alpha-\frac{s_0}2\,\alpha\,(\alpha-1)\,(1-s_0)^{\alpha-2}\quad&\mbox{if}\quad\alpha\leq 2\;,
\end{array}\right.
\]
so that, with $\alpha=1/(1-m)$, $s=-(D_0-D_*)/\big(D_0+\frac{1-m}{2\,m}\,|x|^2\big)\leq -s_0=\frac {D_*}{D_0}-1$, we get the existence of a positive constant $\mathcal M_0=\mathcal M_0(m,D_*,D_0)$ such that
\[
w(t,x)-1\ge\frac{V_{D_0}}{V_{D_*}}-1=\left(1-\frac{D_0-D_*}{D_0+\frac{1-m}{2\,m}\,|x|^2}\right)^\frac 1{1-m}-1\ge\frac{\mathcal M_0}{D_0+\frac{1-m}{2\,m}\,|x|^2}\quad\forall\; x\in\RR^d\;.
\]
Hence, there exists a positive constant $\mathcal M$ depending on $\max\{\mathcal M_0,\mathcal M_1\}$, $D_0$ and $D_1$, for which we obtain
\[
|w(t,x)-1|\le \min\left\{\frac{\sigma_0\,\ee^{-\gamma_\infty\,(t-t_0)}}{V_{D_*}}\,,\;\mathcal M\,V_{D_*}^{1-m}\right\}=\sigma\,\ee^{-\gamma_\infty\,\frac{1-m}{2-m}\,(t-t_0)}\,,\quad\sigma:=\mathcal M^\frac 1{2-m}\,\sigma_0^\frac{1-m}{2-m}\,.
\]
As a consequence, for any $t\geq t_0$, we have improved bounds on $w$, with $W_0$ and $W_1$ replaced respectively by $\sigma_0(t):=1-\sigma\,\ee^{-\gamma_\infty\frac{1-m}{2-m}\,(t-t_0)}$ and $\sigma_1(t):=1+\sigma\,\ee^{-\gamma_\infty\frac{1-m}{2-m}\,(t-t_0)}$. As in Proposition~\ref{Thm:expodecfreeenergy}, according to Theorem~\ref{Cor:E-I} and Inequality~\eqref{Eqn:ExprGamma}, $z(t)=\fw{w(t)}$ satisfies
\[
\frac{dz}{dt}\leq - \gamma(t)\,z(t)
\]
with
\[
\gamma(t):=2\,\frac{m\,\sigma_0(t)^{2(2-m)}-\Cmd\,d\,(1-m)\left[ \sigma_1(t)^{2(2-m)}-\sigma_0(t)^{2(2-m)}\right]}{\Cmd\,\sigma_0(t)^{1-m}\,\sigma_1(t)^{2(2-m)}}=2\,\gopt-\eta(t)\;,
\]
$\gopt=m/\Cmd$ and $\eta(t)\leq \eta_0\,\ee^{-\gamma_\infty\frac{1-m}{2-m}\,(t-t_0)}$ for some $\eta_0>0$. A Gronwall argument then shows that for any $t\geq t_0$,
\[
\log\left(\frac{z(t)}{z(t_0)}\right)\leq-\,2\,\gopt\,(t-t_0)+\frac{\eta_0}{\gamma_\infty\frac{1-m}{2-m}}\,\left[1-\ee^{-{\gamma_\infty\frac{1-m}{2-m}}(t-t_0)}\right]\leq-\,2\,\gopt\,(t-t_0)+\frac{\eta_0\,(2-m)}{\gamma_\infty\,(1-m)}\;.
\]
which completes the estimate on $\fw {w(t)}$. For $t$ large enough, $\frac 12\,\gamma(t)\in(\lambda,\gopt)$ and the $\LL^\infty$ estimate follows.
\qed

\subsubsection*{Proof of Theorem~\ref{thm:CRE-exp}}
As in the proof of Corollary~\ref{cor3.4bis}, a H\"older interpolation inequality shows that, for any $\delta>0$,
\[
\int_{\RR^d}\left|w(t)-1\right|^{p}\dx\le\|w(t)-1\|_\infty^{p-\delta-\frac d2}\;\int_{\RR^d}\left(\frac{c_0}{1+|x|^2}\right)^{\delta+\frac d2}\!\!\dx\;.
\]
\qed

\subsection{Proof of Theorem~\ref{Thm:A2}}

We first apply Theorem~\ref{Thm:6.3} and Lemma~\ref{Conv.Weight}, (i), with $\vartheta=0$ to obtain, for some $t_0\geq 0$,
\begin{equation*}
\left\|v(t)-V_{D_*}\right\|_{2}\leq K_{\vartheta=0}\,\left( \fer {v(t)}{V_{D_*}}\,\right)^{\frac12} \le C_2\,\ee^{-\gopt\,t}\quad\forall\; t\geq t_0\;,
\end{equation*}
for some positive constant $C_2$. By the interpolation inequality~\eqref{eq:interpolation} with $\lambda=-\alpha\,d< 0=\mu<1/2=\nu$, $C=\mathcal C_{-\alpha d,\,0,\,1/2}$, and Lemma~\ref{Lemma:Calpha}, \eqref{eq:mathcalH}, we have
\begin{equation*}
\|v(t)-V_{D_*}\|_\infty \le\,C\,\|v(t)-V_{D_*}\|_{C^\alpha}^\theta\; \|v(t)-V_{D_*}\|_2^{1-\theta} \le\,C\,\mathcal{H}^\theta\,\|v(t)-V_{D_*}\|_\infty^\theta\;\|v(t)-V_{D_*}\|_2^{1-\theta}
\end{equation*}
where $\theta=1/(2+\alpha\,d)$. This implies
\begin{equation*}
\|v(t)-V_{D_*}\|_\infty\le C^{1/(1-\theta)}\,\mathcal{H}^{\theta/(1-\theta)}\,\|v(t)-V_{D_*}\|_2\quad\forall\; t\geq t_0\;.
\end{equation*}
{}From H\"older's inequality, $\|v(t)-V_{D_*}\|_q\le \|v(t)-V_{D_*}\|_\infty^{(q-2)/q}\, \|v(t)-V_{D_*}\|_2^{2/q}$, $q\in(2,\infty]$, we deduce that $\|v(t)-V_{D_*}\|_q$ decays with the same rate as $\|v(t)-V_{D_*}\|_2$. If $q\in(q_*,2)$, we apply Lemma~\ref{Conv.Weight}, (ii), and Theorem~\ref{Thm:6.3} to prove that for some positive constant $C_q$ and for some $t_0\geq 0$,
\begin{equation*}
\|v(t)-V_{D_*}\|_q\le C_q\,\ee^{-\gopt\,t}\quad\forall\; t\geq t_0\;.
\end{equation*}
Similarly, the estimate $\|v(t)-V_{D_*}\|_{C^j(\RR^d)}$ follows from Lemma~\ref{Conv.Weight}, (iii), and Theorem~\ref{Thm:6.3}. This completes the proof of Theorem~\ref{Thm:A2}. \qed

\section*{Appendix A: Hardy-Poincar\'e inequalities}
\setcounter{subsection}{0}
\setcounter{equation}{0}
\setcounter{theorem}{0}
\renewcommand{\thesection}{A}

In this appendix, we state and prove a result on inequalities which we have already been partially studied in~\cite{BBDGV-CRAS}. Here we give more details and a few improvements. We are especially interested in the explicit values of the constants which enter in the convergence rates of Theorems~\ref{Thm:A2} and~\ref{thm:CRE-exp}. This is why we take weights which are adapted to Equation \eqref{FP.Problem} and define the measures
\[
\dmu:=V_D^{2-m}\dx\quad\mbox{and}\quad\dnu:=V_D\dx\;,
\]
where $V_D(x)=\left(D+\frac{1-m}{2\,m}\,|x|^2\right)^{-1/(1-m)}$. Incidentally we observe that $ \dmu=V_D^{1-m}\,\dnu$. To a function $g\in \LL1(\dmu)$, we associate its average $\overline{g}=\int_{\RR^d}g(x)\;\dmu$. Recall that $m_*=(d-4)/(d-2)$.

\subsection{Statement and comments}\label{Sec:A0}

\begin{theorem}\label{Teor.Spect.} Let $d\geq 1$ and $D> 0$. If $m\in(0,1)$ and $1\leq d\leq 4$, or $m\in(m_*,1)$ and $d\geq 5$, then there exists a positive constant $\Cmd$, which does not depend on $D$, such that
\begin{equation}\label{Gap.Supercrit.Thm}
\int_{\RR^d}\left|g-\overline{g}\right|^2\,\dmu\le \Cmd\int_{\RR^d}\left|\nabla g\right|^2\dnu\quad\forall\;g\in\mathcal D(\RR^d)\;,\quad\overline{g}=\int_{\RR^d}g\;d\mu\;.
\end{equation}
In case $d\geq 5$ and $m\in(0,m_*)$, we have
\begin{equation}\label{Gap.Subcrit.Thm}
\int_{\RR^d}g2\,\dmu\le \Cmd\int_{\RR^d}\left|\nabla g\right|^2\,\dnu\quad\forall\;g\in\mathcal D(\RR^d)
\end{equation}
and $\Cmd=\frac{8\,m\,(1-m)}{[(d-2)\,(m-m_*)]2}$ is optimal. \end{theorem}
Estimates of the optimal constant $\Cmd$ when $m>m_*$ are given
below in Proposition~\ref{Thm3b}. With
$v_m(x)=(1+|x|^2)^{-1/(1-m)}$, a simple change of variables shows
that $\gopt=m/\Cmd$ is such that
 \begin{equation}\label{DeMuav_m}
\gopt\;=\;m\;\inf_h\;\frac{\int_{\RR^d}|\nabla
h|^2\,V_{D}\dx}{\int_{\RR^d}|h-\bar
h|^2\,V_{D}^{2-m}\dx}\;=\;\frac{1-m}2\,\inf_h\frac{\int_{\RR^d}|\nabla
h|^2\,v_m\dx}{\int_{\RR^d}|h-\tilde h|^2\,v_m^{2-m}\dx}\;,
\end{equation}
where the infima are taken over the set of smooth functions $h$ such that\\
- either $m<m_*$ and ${\rm supp}(h)\subset\RR^d\setminus\{0\}$ and $\bar h=0$, $\tilde h=0$, \\
- or $m>m_*$,
\[
\bar
h:=\frac{\int_{\RR^d}h\,V_{D}^{2-m}\dx}{\int_{\RR^d}V_{D}^{2-m}\dx}\quad\mbox{and}\quad\tilde
h:=\frac{\int_{\RR^d}h\,v_m^{2-m}\dx}{\int_{\RR^d}v_m^{2-m}\dx}\;.
\]
This already shows that $\gopt$ is independent of $D$.

We observe that as $|x|\to\infty$, $\dmu\sim\dnu/|x|^2$. Hence, if
$m\in(0,m_*)$, Inequality~\eqref{Gap.Subcrit.Thm} is of Hardy
type. Otherwise, if $m\in(m_*,1)$,
Inequality~\eqref{Gap.Supercrit.Thm} involves an average and is
rather of Poincar\'e type. In such a case, we shall also say that
it is a weighted Poincar\'e inequality, or that there is a
spectral gap, since for the associated operator, the lowest
eigenvalue, $0$, is achieved by the constant functions, and the
second eigenvalue corresponds to $\gopt=m/\Cmd$ where $\mathcal
C_{m,d}$ is the best constant in the inequality.
See~\cite{BBDGV-CRAS} for further considerations on these issues.

We also remark that Theorem~\ref{Teor.Spect.} provides an explicit
example for which the weighted Poincar\'e inequality holds, while
the corresponding weighted logarithmic Sobolev inequality does not
hold, even in dimension $d=1$, as shown by~\cite[Theorem~3]{BR1}.

The proofs of~\eqref{Gap.Supercrit.Thm}
and~\eqref{Gap.Subcrit.Thm} are quite different and for this
reason we treat the two cases separately. We start with the proof
of \eqref{Gap.Subcrit.Thm} corresponding to the case $m<m_*$,
$d\geq 5$.

\subsection{Case $m\in(0,m_*)$}\label{Sec:A1}

The proof follows the ideas of~\cite{BBDGV-CRAS}. We reproduce it here for completeness. We compute
\[
|\nabla V_D(x)|^2=\frac{|x|^2}{m2}\,V_D(x)^{2(2-m)}
\]
and
\[
-\,2\,m2\,\frac{\Delta V_D(x)}{V_D(x)^{3-2m}}=2\,d\,D\,m+(d-2)\,(m_*-m)|x|^2\;.
\]
An integration by parts and the Cauchy-Schwarz inequality show that

\[\begin{split}
\left|\int_{\RR^d}|g|^2\,\Delta V_D\,\rd x\right| &\le 2\int_{\RR^d}|g|\,|\nabla g|\,|\nabla V_D|\,\rd x\\
&\le 2\left(\int_{\RR^d}|g|^2\,|\Delta V_D|\,\rd x\right)^{1/2}\left(\int_{\RR^d}|\nabla g|^2\,|\nabla V_D|^2\,|\Delta V_D|^{-1}\,\rd x\right)^{1/2}.
\end{split}
\]
As in~\cite{MR1612685}, we remark that $\Delta V_D$ has a constant sign and get the estimate
\[
\left|\int_{\RR^d}|g|^2\,\Delta V_D\,\rd x\right|=\int_{\RR^d}|g|^2\,|\Delta V_D|\,\rd x\le 4\int_{\RR^d}|\nabla g|^2\,|\nabla V_D|^2\,|\Delta V_D|^{-1}\,\rd x\;.
\]
Weights can be estimated on both sides of the inequality:
\begin{eqnarray*}
&& \frac{|\Delta V_D|}{V_D^{2-m}}=\frac{2\,d\,D\,m+(d-2)\,(m_*-m)\,|x|^2}{m\,(2\,D\,m+(1-m)\,|x|^2)} \ge \frac{(d-2)\,(m_*-m)}{m\,(1-m)}\;,\\
&&\frac{|\nabla V_D|^2}{|\Delta V_D|\,V_D}\le\frac{2\,|x|^2}{2\,d\,D\,m+(d-2)\,(m_*-m)\,|x|^2}\le\frac{2}{(d-2)\,(m_*-m)}\;,
\end{eqnarray*}
which proves \eqref{Gap.Subcrit.Thm}. See~\cite{BBDGV-CRAS} for further details.

\medskip We now consider the limit $D\to 0^+$. With $\alpha:=1/(m-1)\in(1-d/2,-1)$, that is $m\in(0,m_*)$, and
\[
\kappa_\alpha:=\frac{8\,m\,(1-m)}{[(d-2)\,(m-m_*)]2}\cdot\frac{1-m}{2\,m}=\frac{4\,(1-m)2}{[(d-4)-(d-2)\,m]2}\;,
\]
Inequality~\eqref{Gap.Subcrit.Thm} takes the form of a weighted Hardy inequality,
\begin{equation*}\label{Hardy.Weighted.Alpha0}
\int_{\RR^d}\frac{|g|^2}{|x|^2}\;|x|^\alpha\dx\le\kappa_\alpha\int_{\RR^d}|\nabla g|^2\,|x|^\alpha\dx\quad\forall\;g\in\mathcal D(\RR^d)\;.
\end{equation*}
Such an inequality is easy to establish by the ``completing the square method'' as follows. Let $\alpha\in\RR\setminus\{\alpha\}$ with $\alpha_*:=1-d/2$, and $g\in \mathcal D(\RR^d)$. Then
\begin{eqnarray*}
0&\leq &\int_{\RR^d}\left|\,\nabla g+\lambda\,\frac x{|x|^2}\,g\,\right|^2\,|x|^{2\alpha}\dx\\
&&=\int_{\RR^d}|\nabla g|^2\,|x|^{2\alpha}\dx+\Big[\lambda2-\lambda\,(2\,\alpha+d-2)\Big]\int_{\RR^d}\frac{|g|^2}{|x|^2}\;|x|^{2\alpha}\dx\;.
\end{eqnarray*}
An optimization of the right hand side with respect to $\lambda$ results in choosing $\lambda=(2\,\alpha+d-2)/2$, that is
\[
\frac{1}{\lambda2}=\frac{4}{(d+2\,\alpha-2)2}=\frac{4\,(1-m)2}{[(d-4)-(d-2)\,m]2}=\kappa_\alpha\;.
\]
The weighted Hardy inequality is optimal, with optimal constant
$\kappa_\alpha$, as follows by considering the test functions
$g_\varepsilon(x):=\min\{\varepsilon^{-\lambda},(|x|^{-\lambda}-\varepsilon^{\lambda})_+\}$
and letting $\varepsilon\to0$.\qed

A closer inspection of the proof reveals that the constant $\kappa_\alpha$ in the weighted
Hardy inequality also is optimal when $m>m_*$. Consider indeed the test functions
$g_\varepsilon(x):=|x|^{1-\alpha-d/2+\varepsilon}$ for $|x|<1$ and $g_\varepsilon(x)=(2-|x|)_+$
for $|x|\geq 1$, and then let $\varepsilon\to0$.
\begin{proposition}[Weighted Hardy inequality]\label{Cor:WeightedHardy}
With the above notations, for any $\alpha\in\RR$, $\alpha\ne\alpha_*$,
\begin{equation*}\label{Hardy.Weighted.Alpha}
\int_{\RR^d}\frac{|g|^2}{|x|^2}\;|x|^{2\alpha}\dx\le\kappa_\alpha\int_{\RR^d}|\nabla g|^2\,|x|^{2\alpha}\dx\quad\forall\;g\in\mathcal D(\RR^d)\;,
\end{equation*}
with the additional requirement that $g$ is supported in
$\RR^d\setminus\{0\}$ if $\alpha<\alpha_*$, and $\kappa_\alpha$ is
optimal. \end{proposition}
The range $m\in(0,1)$ corresponds to
$1/(m-1)=\alpha\in(-\infty,-1)$, so that $m=m_*$ is equivalent to
$\alpha=\alpha_*$. Notice that the result holds without other
restriction than $\alpha\neq\alpha_*$, but one has to be careful
with integrability condition at $x=0$ if $\alpha<\alpha_*$.

\subsection{Case $\max\{0,m_*\}<m<1$}\label{Sec:A2}

Several partial results are known. In the range $m\in(m_c,1)$, see
\cite{MR1901093} for an estimate of $\Cmd$ based on the
Bakry-Emery method, \cite{BBDGV-CRAS} for other estimates,
and~\cite{MR2126633} for the exact values of the optimal constant
for a  corresponding  linear problem.

We now prove~\eqref{Gap.Supercrit.Thm} with some explicit
estimates of the constant $\Cmd$ in the whole range
$(\max\{0,m_*\},(d-2)/(d-1))\supset(m_*,m_c]$. Because of the
change of variables \eqref{DeMuav_m}, our task is now to
characterize $\Cmd$ as
\[
\left(\frac{(1-m)\,\Cmd}{2\,m}\right)^{-1}\;=\;\inf_h\frac{\int_{\RR^d}|\nabla h|^2\,v_m\dx}{\int_{\RR^d}|h-\tilde h|^2\,v_m^{2-m}\dx}\;.
\]On~$\RR^+$, consider the function $\m(r):=r^{d-1}(1+|r|^2)^{(2-m)/(m-1)}$, and denote its median by $\eta$. Let $\nu(r):=r^{d-1}(1+r2)^{1/(m-1)}$ and define for all $\zeta>0$ the quantity
\begin{equation*}
\textsf K(\zeta):=\frac{2\,m}{1-m}\,\max\left\{\mbox{\textsf A}(\zeta),\mbox{\textsf B}(\zeta)\right\}
\end{equation*}
with
\[
\mbox{\textsf A}(\zeta):=\sup_{r<\zeta}\left[\int_0^{r}\m(s)\ds\int_r^\zeta\frac{\ds}{\nu(s)}\right]\,,\quad
\mbox{\textsf B}(\zeta):=\sup_{r>\zeta} \left[\int_\zeta^r\frac{\ds}{\nu(s)}\int_r^{+\infty}\mu(s)\ds\right]\,.
\]
By convention, we take $\mbox{\textsf K}(0)=\frac{2\,m}{1-m}\,\mbox{\textsf B}(0)$. The following result is inspired by~\cite{BR1, bobkov-gotze, MR1777034, Miclo}.
\begin{proposition}\label{Thm3b} Let $d\geq 1$. For any $m\neq m_*$,
\[
\Cmd\geq \frac{8\,m\,(1-m)}{\left[d-4-m\,(d-2)\right]2}\;.
\]
If $m\in(m_*,1)$, then
\[\label{Const.Gap.Supercrit.Thm}
\mathcal{C}_{m,1}\leq\mbox{\rm\textsf K}(0)\quad\mbox{and}\quad\Cmd\leq\max\left\{2\,\mbox{\rm\textsf K}(\eta)\,,\;\frac{4\,m}{(1-m)\,(d-1)}\right\}\quad\mbox{if}\quad d\geq 2\;,
\]
where, for any $m\in(m_*,(d-2)/(d-1))$,
\begin{equation*}
\mbox{\rm\textsf K}(\eta)\le \frac{m\,(2-m)\,2^{\frac{3-2m}{1-m}}\left(1+2^{\frac{2-m}{1-m}}\right)}{d\,\left[d-4-m\,(d-2)\right]2}\;.
\]\end{proposition}
The function $v_m^{2-m}$ is integrable for any $m\in (m_*,1)$, so that $\mbox{\rm\textsf K}(\eta)$ is well defined in this range. The upper bound on $\Cmd$ is equal to its exact value up to a factor which is at least $1/4$ (and at most $1$). Such an interval is inherent to the method, see~\cite{Miclo}. The case $m=m_c\leq (d-2)/(d-1)$ is covered, showing in particular that $\mathcal C_{m_c,d}$ is positive, finite. The bounds diverge as $m\searrow m^*$ with same behavior at first order. Our approach can be extended easily to the case $((d-2)/(d-1),1)$, with slightly different estimates of $\textsf K(\eta)$, but this case is already covered in \cite{BBDGV-CRAS,MR1901093,MR2126633} by other methods. The restriction $m<(d-2)/(d-1)$ is convenient from a technical point of view, and not essential at all.
\proof The lower bound on $\Cmd$ is achieved as in Section~\ref{Sec:A1} by taking the limit $D=0$, thus showing that $\Cmd\geq \frac{1-m}{2\,m}\,\kappa_{1/(m-1)}$ and using Proposition~\ref{Cor:WeightedHardy}.

Let us prove the upper bounds. We introduce the standard change of variables from Cartesian to spherical coordinates, \ie\ $r=|x|$, and $\vartheta=x/|x|$. In these coordinates, the gradient can be written as $(\partial_r,\frac 1r\nabla_\theta)$ where $\partial_r=\frac xr\cdot\nabla$ is the partial derivative with respect to the radial variable $r$ and $\nabla_\theta$ is the derivative with respect to the angular variables. We shall denote by $\mathbb S^{d-1}\subset\RR^d$ the unit sphere and parametrize it with the variable $\vartheta$.

The radial density functions $r\mapsto\mu(r)$ and $r\mapsto\nu(r)$ are such that $v_m\!\dx=\mu(|x|)\,\rd x$ and $v_m^{2-m}\dx=\nu(|x|)\,\rd x$. We introduce the following normalization constants:
\[
\omega_d=\int_{\mathbb{S}^{d-1}}\rd\vartheta=\frac{2\,\pi^{d/2}}{\Gamma(d/2)}\;,\qquad\widehat{\rd\vartheta}=\omega_d^{-1}\,\rd\vartheta\;,\qquad \int_{\mathbb{S}^{d-1}}\widehat{\rd\vartheta}=1\;.
\]
With these notations,
\[
\mu(r)\,\rd r\,\rd\vartheta=\frac{v_m(x)}{1+|x|^2}\dx=v_m^{2-m}\dx\quad\mbox{and}\quad\nu(r)\,\rd r\,\rd\vartheta=v_m\dx\;.
\]
We define a directional average of a function $f$ by
\[
\overlinetilde{f_\mu}(\vartheta):=\int_0^{+\infty}f(r,\vartheta)\,\hat\mu(r)\,\rd
r\quad\mbox{with}\quad\hat\mu(r):=\frac{\mu(r)}{\int_0^{+\infty}\mu(s)\dd
s}
\]
and the global average of $f$ by
\[
\overlinetilde{f}:=\frac{\int_{\RR^d}f\,v_m^{2-m}\dx}{\int_{\RR^d}v_m^{2-m}\dx}
=\iintpolar f(r,\vartheta)\,\hat\mu(r)\,\rd
r\,\widehat{\rd\vartheta}
=\int_{\mathbb{S}^{d-1}}\overlinetilde{f_\mu}(\vartheta)\,\widehat{\rd\vartheta}\;.
\]

In the case $d=1$, Theorem 2 of~\cite{BR2}, also see~\cite{Miclo},
says that Inequality~\eqref{Gap.Supercrit.Thm} holds with
\[
\mathcal{C}_{m,1}\leq\frac{2\,m}{1-m}\,\sup_{r>0} \left[\int_r^{+\infty}\mu(r)\dd r\,\int_0^r\frac{\rd r}{\nu(r)}\right]=\mbox{\rm\textsf K}(0)
\]
in which case we also have the estimate
$\mathcal{C}_{m,1}\le\mbox{\rm\textsf K}(0)\le
4\,\mathcal{C}_{m,1}$.

In case of radial functions, Inequality~\eqref{Gap.Supercrit.Thm} takes the form:
\begin{equation}\label{Gap.1dim}
\int_0^{+\infty}\left|f(r)2-\overlinetilde{f}\,\right|^2 \,\mu(r)\,\rd r\le \frac{1-m}{2\,m}\,\Cmd^{\rm rad}\int_0^{+\infty}\left|f^\prime(r)\right|^2\,\nu(r)\,\rd r
\end{equation}
with
\begin{equation*}\label{Const.Gap.1dim}
\Cmd^{\rm rad}\leq \frac{2\,m}{1-m}\,\max\left\{\sup_{r>\eta} \left[\int_r^{+\infty}\mu(r)\dd r\,\int_\eta^r\frac{\rd r}{\nu(r)}\right]\;,\;\sup_{r<\eta}\left[\int_0^r\mu(r)\dd r\int_r^\eta\frac{\rd r}{\nu(r)}\right] \right\}=\mbox{\rm\textsf K}(\eta)\;.
\end{equation*}
It is straightforward to show that $\mbox{\rm\textsf K}(\eta)$ is finite, with the present choices of $\mu$ and $\nu$, for $m\in(m_*,(d-2)/(d-1))$, and as above, $\Cmd^{\rm rad}\le\mbox{\rm\textsf K}(\eta)\le4\,\Cmd^{\rm rad}$.

\medskip We now focus on the case of non radial functions, with $d\ge 2$, and rewrite the left hand side of~\eqref{Gap.Supercrit.Thm} in spherical coordinates.
\begin{eqnarray*}\label{n.dim.Gap}
\int_{\RR^d}\left|f(x)-\overlinetilde{f}\,\right|^2\,v_m^{2-m}\dx&=&\omega_d\,\iintpolar\big|f(r,\vartheta)-\overlinetilde{f}\,\big|^2\,\mu(r)\,\rd r\,\widehat{\rd\vartheta}\\
&=&\omega_d\,\iintpolar\big|f(r,\vartheta)-\overlinetilde{f_\mu}(\vartheta) +\overlinetilde{f_\mu}(\vartheta)-\overlinetilde{f}\,\big|^2\,\mu(r)\,\rd r\,\widehat{\rd\vartheta}\\
&\le& 2\,\omega_d\,\big[\mbox{\rm (I)}+\mbox{\rm (II)}\big]
\end{eqnarray*}
with
\begin{eqnarray*}
\mbox{\rm (I)}&=&\iintpolar\big|f(r,\vartheta)-\overlinetilde{f_\mu}(\vartheta)\big|^2\,\mu(r)\,\rd r\,\widehat{\rd\vartheta}\;,\\
\mbox{\rm (II)}&=&\iintpolar\big|\overlinetilde{f_\mu}-\overlinetilde{f}\,\big|^2\,\mu(r)\,\rd r\,\widehat{\rd\vartheta}=\int_{\mathbb{S}^{d-1}}\left|\overlinetilde{f_\mu}-\overlinetilde{f}\,\right|^2\widehat{\rd\vartheta}\;.
\end{eqnarray*}
We estimate (I) by~\eqref{Gap.1dim} and get
\begin{equation*}\label{radial.Gap}
\mbox{\rm (I)}=\iintpolar\big|f(r)2-\overlinetilde{f_\mu}(\vartheta)\big|^2 \,\mu(r)\,\rd r\,\widehat{\rd\vartheta} \le \frac{1-m}{2\,m}\,\Cmd^{\rm rad}\iintpolar\big|\partial_r\,f(r,\vartheta)\big|^2\,\nu(r)\,\rd r\,\widehat{\rd\vartheta}\;.
\end{equation*}
To estimate (II), we rely on the Poincar\'e inequality on the unit sphere $\mathbb{S}^{d-1}$,
\begin{equation*}\label{Poinc.Sphere}
\int_{\mathbb{S}^{d-1}}\left|u-\hat u\,\right|^2\,\widehat{\rd\vartheta}\le \frac 1{d-1}\int_{\mathbb{S}^{d-1}}\left|\nabla_\vartheta\,u\right|^2\,\widehat{\rd\vartheta}\quad\forall\;u\in H1(\mathbb{S}^{d-1})\;.
\end{equation*}
Here $\hat u:=\int_{\mathbb{S}^{d-1}}u\,\widehat{\rd\vartheta}$. In the inequality, $1/(d-1)$ is the optimal constant, as can be checked using spherical harmonic functions. See for instance~\cite{MR1164616,BL,SC1}. The inequality itself can be recovered by various methods. For example, using the inverse stereographic projection, see~\cite{MR717827}, the optimal Sobolev inequality on $\RR^d$ becomes
\[
\left(\int_{\mathbb{S}^{d-1}}|v|^p\,\widehat{\rd\vartheta}\right)^{2/p}\le\int_{\mathbb{S}^{d-1}}|v|^2\,\widehat{\rd\vartheta}+\frac{p-2}{d-1}\int_{\mathbb{S}^{d-1}}\left|\nabla_\vartheta\,v\right|^2\,\widehat{\rd\vartheta}\;,
\]
for any $u\in H1(\mathbb{S}^{d-1})$, with $p=2\,d/(d-2)$, $d\geq 3$. The inequality also holds true for any $p\in(2,2\,d/(d-2))$ if $d\geq 3$ and for any $p>2$ if $d=2$, see~\cite{MR1230930}. Hence we recover the Poincar\'e inequality on $\mathbb{S}^{d-1}$ by writing $v=1+\varepsilon\,u$ and keeping only the terms of order $\varepsilon2$ as $\varepsilon\to 0$.

We apply the Poincar\'e inequality with $u=\overlinetilde{f_{\mu}}$.
\begin{equation*}
\int_{\mathbb{S}^{d-1}}\left|\overlinetilde{f_\mu}-\overlinetilde{f}\,\right|^2\widehat{\rd\vartheta}\le\frac 1{d-1}\int_{\mathbb{S}^{d-1}}\left|\nabla_\vartheta\,\overlinetilde{f_\mu}\right|^2\widehat{\rd\vartheta}
\end{equation*}
Recall that $|\nabla f|^2= |\partial_r f(r,\vartheta)|^2+\frac{1}{r2}\,|{\nabla_\vartheta}\,f(r,\vartheta)|^2$. Using the Cauchy-Schwarz inequality and the estimate $r2\,\mu(r)\,\rd r\leq\nu(r)\,\rd r$, we get
\begin{eqnarray*}
\int_{\mathbb{S}^{d-1}}\left|\nabla_\vartheta\,\overlinetilde{f_\mu}\right|^2\widehat{\rd\vartheta} &\leq& \iintpolar\big|\nabla_\vartheta\,f(r,\vartheta)\big|^2\,\mu(r)\,\rd r\,\widehat{\rd\vartheta}\\
&\leq&\iintpolar\kern 25pt\frac{1}{r2}\left|\nabla_\vartheta\,f(r,\vartheta)\right|^2\,\nu(r)\,\rd r\,\widehat{\rd\vartheta}\;.
\end{eqnarray*}
This proves that
\[
\mbox{\rm (II)}\leq \frac 1{d-1}\,\iintpolar\kern 25pt\frac{1}{r2}\left|\nabla_\vartheta\,f(r,\vartheta)\right|^2\,\nu(r)\,\rd r\,\widehat{\rd\vartheta}\;.
\]

Summarizing, we have shown that
\[
\int_{\RR^d}\left|f(x)-\overlinetilde{f}\,\right|^2\,v_m^{2-m}\dx\leq 2\,\max\left\{\frac{1-m}{2\,m}\,\Cmd^{\rm rad},\frac 1{d-1}\right\}\int_{\RR^d}|\nabla f|^2\,v_m\dx\;.
\]
By undoing the change of coordinates as in \eqref{DeMuav_m}, we get
\[
\int_{\RR^d}\left|f(x)-\overline{f}\right|^2\rd\mu \leq \max\left\{2\,\Cmd^{\rm rad},\frac{4\,m}{(1-m)\,(d-1)}\right\}\int_{\RR^d}|\nabla f|^2\,\rd\nu\;.
\]
The bounds on ${\textsf K}(\eta)$ follow by quite long but
straightforward calculations, omitted here. \qed

\section*{Appendix B: Extension to exponents $m\le 0$}
\setcounter{subsection}{0}
\setcounter{equation}{0}
\renewcommand{\thesection}{B}

The presence at several instances of factors of the form $1/m$ in
the previous calculations may suggest that there is an essential
divergence as $m\to 0$. In the present section we want to dispel
that impression by introducing a normalization that is often used
in the literature, consisting in a rescaling of the time variable
of the form $\tau'=m\,\tau$ that modifies equation
\eqref{FDE.Problem} into
\begin{equation}\label{MFDE}
\partial_{\tau'} u=\nabla \cdot (u^{m-1}\nabla u)\;.
\end{equation}
One of the first consequences is that the new equation, that we
will call {\sl modified fast diffusion equation\/} for clarity,
makes perfect sense as a nonlinear parabolic equation of singular
type for all the range of exponents $m\in \RR$ (including $m=0$),
in particular for all $m<1$ that form the extended range of the
fast diffusion. Such approach has been consistently used in
\cite{Vazquez2006} where it is shown that the effect on the
self-similar solutions of Barenblatt type  is just to eliminate
the denominator $m$ in the formulas \eqref{baren.form1},
\eqref{baren.form2}. Note the rescaling to obtain \eqref{MFDE}
from the standard fast diffusion equation \eqref{FDE.Problem}
 when $m>0$ can also be done by changing the space
variable in the form $x=\sqrt{m}\,x'$ and not changing time.

Since the general theory (existence, uniqueness, estimates,
special solutions and extinction) has  been developed to the
measure we need it, we can follow the different stages of the
present paper with due attention to chasing the $m$ factors, and
the results stated in Section~\ref{Sec:Intro} remain valid. For
instance, the formula defining $\fw w$ in Section \ref{rel.entr}
has to be replaced by
\begin{equation*}
\fw w=-\int_{\RR^d}\big[\,\log w-(w-1)\big]\dx\;.
\end{equation*}
In the linearization of Section \ref{Sec:SpectralGap} there is no
$m$ factor in the definition of operator $A_m$ and neither in the
definition of $\IL g$. Let us  mention two other points of
interest: the exponent $m_*$ becomes negative for $d=1,2,3$ and
zero for $d=4$, but it still plays the same role of an important
critical exponent separating different behavior types. On the
other hand, the constant $\gopt$ that gives the decay rate in our
main result has a finite positive value as $m\to 0$, according to
Theorem~\ref{Teor.Spect.}.  This constant determines the rates in
all results concerning the asymptotic behavior of the solutions
 and is not affected by our $m$-rescaling.

We have refrained from treating the extension to $m\le 0$ in the
main body of the paper in order to avoid further distractions in
an already very technical matter. Whole details will
appear separately.

\bigskip\noindent{\small{\sc Acknowledgments.} {This project has been supported by the IFO project of the French Research Agency (ANR). A.B. acknowledges the support of a bourse Lavoisier. M.B. thanks CNRS, CEREMADE, Dpto. di Matematica of Politecnico di Torino and the Dpto. de Matem\'{a}ticas of Universidad Aut\'{o}noma de Madrid for post-doctoral grants during which this project has been carried out. J.L.V. was partially supported by Spanish Project MTM2005-08760-C02-01. A.B. and J.L.V. were partially supported by the ESF Programme \lq\lq Global and geometric aspects of nonlinear partial differential equations".}}

\par\smallskip\noindent{\small \copyright~2007 by the authors. This paper may be reproduced, in its entirety, for non-commercial~purposes.}
\def\cprime{$'$}


\begin{thebibliography}{10}

\bibitem{MR0448404}
{\sc T.~Aubin}, {\em Probl\`emes isop\'erim\'etriques et espaces de {S}obolev},
  J. Differential Geometry, 11 (1976), pp.~573--598.

\bibitem{MR889476}
{\sc D.~Bakry and M.~{\'E}mery}, {\em Diffusions hypercontractives}, in
  S\'eminaire de probabilit\'es, XIX, 1983/84, vol.~1123 of Lecture Notes in
  Math., Springer, Berlin, 1985, pp.~177--206.

\bibitem{BR1}
{\sc F.~Barthe and C.~Roberto}, {\em Sobolev inequalities for probability
  measures on the real line}, Studia Math., 159 (2003).

\bibitem{BR2}
{\sc F.~Barthe and C.~Roberto}, {\em Modified logarithmic {S}obolev
  inequalities on $\mathbb{R}$}, 2006.

\bibitem{MR1164616}
{\sc W.~Beckner}, {\em Sobolev inequalities, the {P}oisson semigroup, and
  analysis on the sphere {$S\sp n$}}, Proc. Nat. Acad. Sci. U.S.A., 89 (1992),
  pp.~4816--4819.

\bibitem{MR1230930}
\leavevmode\vrule height 2pt depth -1.6pt width 23pt, {\em Sharp {S}obolev
  inequalities on the sphere and the {M}oser-{T}rudinger inequality}, Ann. of
  Math. (2), 138 (1993), pp.~213--242.

\bibitem{BBDGV-CRAS}
{\sc A.~Blanchet, M.~Bonforte, J.~Dolbeault, G.~Grillo, and J.-L. V{\'a}zquez},
  {\em {H}ardy-{P}oincar{\'e} inequalities and applications to nonlinear
  diffusions}, C. R. Math. Acad. Sci. Paris,  (2007).

\bibitem{bobkov-gotze}
{\sc S.~G. Bobkov and F.~G{\"o}tze}, {\em Exponential integrability and
  transportation cost related to logarithmic {S}obolev inequalities}, J. Funct.
  Anal., 163 (1999), pp.~1--28.

\bibitem{BL}
{\sc L.~Bo{\v{c}}ek}, {\em Eine {V}ersch\"arfung der
  {P}oincar\'e-{U}ngleichung}, \v Casopis P\v est. Mat., 108 (1983),
  pp.~78--81.

\bibitem{BV-harnackExtTime}
{\sc M.~Bonforte and J.~L. V{\'a}zquez}, {\em Fine asymptotics near extinction
  and elliptic {H}arnack inequalities for the fast diffusion equation}.
\newblock Preprint, April 2006.

\bibitem{BV}
\leavevmode\vrule height 2pt depth -1.6pt width 23pt, {\em Global positivity
  estimates and {H}arnack inequalities for the fast diffusion equation}, J.
  Funct. Anal., 240 (2006), pp.~399--428.

\bibitem{CKN}
{\sc L.~Caffarelli, R.~Kohn, and L.~Nirenberg}, {\em First order interpolation
  inequalities with weights}, Compositio Math., 53 (1984), pp.~259--275.

\bibitem{MR1853037}
{\sc J.~A. Carrillo, A.~J{\"u}ngel, P.~A. Markowich, G.~Toscani, and
  A.~Unterreiter}, {\em Entropy dissipation methods for degenerate parabolic
  problems and generalized {S}obolev inequalities}, Monatsh. Math., 133 (2001),
  pp.~1--82.

\bibitem{MR1901093}
{\sc J.~A. Carrillo, C.~Lederman, P.~A. Markowich, and G.~Toscani}, {\em
  Poincar\'e inequalities for linearizations of very fast diffusion equations},
  Nonlinearity, 15 (2002), pp.~565--580.

\bibitem{MR1777035}
{\sc J.~A. Carrillo and G.~Toscani}, {\em Asymptotic {$L\sp 1$}-decay of
  solutions of the porous medium equation to self-similarity}, Indiana Univ.
  Math. J., 49 (2000), pp.~113--142.

\bibitem{MR1986060}
{\sc J.~A. Carrillo and J.~L. V{\'a}zquez}, {\em Fine asymptotics for fast
  diffusion equations}, Comm. Partial Differential Equations, 28 (2003),
  pp.~1023--1056.

\bibitem{CW}
{\sc F.~Catrina and Z.-Q. Wang}, {\em On the {C}affarelli-{K}ohn-{N}irenberg
  inequalities: sharp constants, existence (and nonexistence), and symmetry of
  extremal functions}, Comm. Pure Appl. Math., 54 (2001), pp.~229--258.

\bibitem{MR1158938}
{\sc Y.~Z. Chen and E.~DiBenedetto}, {\em H\"older estimates of solutions of
  singular parabolic equations with measurable coefficients}, Arch. Rational
  Mech. Anal., 118 (1992), pp.~257--271.

\bibitem{MR1777034}
{\sc S.-K. Chua and R.~L. Wheeden}, {\em Sharp conditions for weighted
  1-dimensional {P}oincar\'e inequalities}, Indiana Univ. Math. J., 49 (2000),
  pp.~143--175.

\bibitem{MR1679782}
{\sc P.~Daskalopoulos and M.~del Pino}, {\em On the {C}auchy problem for {$u\sb
  t=\Delta\log u$} in higher dimensions}, Math. Ann., 313 (1999), pp.~189--206.

\bibitem{Daskalopoulos-Sesum}
{\sc P.~Daskalopoulos and N.~Sesum}, {\em Eternal solutions to the {R}icci flow
  on {$\mathbb R\sp 2$}}, Int. Math. Res. Not.,  (2006), pp.~Art. ID 83610, 20.

\bibitem{Daskalopoulos-Sesum2006}
{\sc P.~Daskalopoulos and N.~Sesum}, {\em On the extinction profile of
  solutions to fast-diffusion}.
\newblock Preprint, 2006.

\bibitem{MR990239}
{\sc E.~B. Davies}, {\em Heat kernels and spectral theory}, vol.~92 of
  Cambridge Tracts in Mathematics, Cambridge University Press, Cambridge, 1989.

\bibitem{MR1612685}
{\sc E.~B. Davies and A.~M. Hinz}, {\em Explicit constants for {R}ellich
  inequalities in {$L\sb p(\Omega)$}}, Math. Z., 227 (1998), pp.~511--523.

\bibitem{MR1940370}
{\sc M.~Del~Pino and J.~Dolbeault}, {\em Best constants for
  {G}agliardo-{N}irenberg inequalities and applications to nonlinear
  diffusions}, J. Math. Pures Appl. (9), 81 (2002), pp.~847--875.

\bibitem{MR1857048}
{\sc M.~Del~Pino and M.~S\'aez}, {\em On the extinction profile for solutions of $u_t=\Delta u^{(N-2)/(N+2)}$}, Indiana Univ. Math. J., 50 (2001), pp.~611--628.

\bibitem{MR2126633}
{\sc J.~Denzler and R.~J. McCann}, {\em Fast diffusion to self-similarity:
  complete spectrum, long-time asymptotics, and numerology}, Arch. Ration.
  Mech. Anal., 175 (2005), pp.~301--342.

\bibitem{fk80}
{\sc A.~Friedman and S.~Kamin}, {\em The asymptotic behavior of gas in an
  {$n$}-dimensional porous medium}, Trans. Amer. Math. Soc., 262 (1980),
  pp.~551--563.

\bibitem{Gr-LSI}
{\sc L.~Gross}, {\em Logarithmic {S}obolev inequalities}, Amer. J. Math., 97
  (1975), pp.~1061--1083.

\bibitem{HP}
{\sc M.~A. Herrero and M.~Pierre}, {\em The {C}auchy problem for {$u\sb
  t=\Delta u\sp m$} when {$0<m<1$}}, Trans. Amer. Math. Soc., 291 (1985),
  pp.~145--158.

\bibitem{MR2162628}
{\sc S.-Y. Hsu}, {\em Classification of radially symmetric self-similar
  solutions of {$u\sb t=\Delta\log u$} in higher dimensions}, Differential
  Integral Equations, 18 (2005), pp.~1175--1192.

\bibitem{MR2194833}
\leavevmode\vrule height 2pt depth -1.6pt width 23pt, {\em Extinction profile
  of solutions of a singular diffusion equation}, Commun. Appl. Anal., 9
  (2005), pp.~67--93.

\bibitem{MR2145602}
\leavevmode\vrule height 2pt depth -1.6pt width 23pt, {\em Large time behavior
  of solutions of a singular diffusion equation in {$\mathbb R\sp n$}},
  Nonlinear Anal., 62 (2005), pp.~195--206.

\bibitem{King}
{\sc J.~King}, {\em {Self-similar behavior for the equation of fast nonlinear
  diffusion.}}, Philos. Trans. R. Soc. Lond., Ser. A, 343 (1993), pp.~337--375.

\bibitem{LSU}
{\sc O.~A. Lady{\v{z}}enskaja, V.~A. Solonnikov, and N.~N. Ural{\cprime}ceva},
  {\em Linear and quasilinear equations of parabolic type}, Translated from the
  Russian by S. Smith. Translations of Mathematical Monographs, Vol. 23,
  American Mathematical Society, Providence, R.I., 1967.

\bibitem{MR1974458}
{\sc C.~Lederman and P.~A. Markowich}, {\em On fast-diffusion equations with
  infinite equilibrium entropy and finite equilibrium mass}, Comm. Partial
  Differential Equations, 28 (2003), pp.~301--332.

\bibitem{MR717827}
{\sc E.~H. Lieb}, {\em Sharp constants in the {H}ardy-{L}ittlewood-{S}obolev
  and related inequalities}, Ann. of Math. (2), 118 (1983), pp.~349--374.

\bibitem{Miclo}
{\sc L.~Miclo}, {\em Quand est-ce que les bornes de {H}ardy permettent de
  calculer une constante de {P}oincar{\'e} exacte sur la droite ?}
\newblock Preprint, 2006.

\bibitem{MR0109940}
{\sc L.~Nirenberg}, {\em On elliptic partial differential equations}, Ann.
  Scuola Norm. Sup. Pisa (3), 13 (1959), pp.~115--162.

\bibitem{MR1842429}
{\sc F.~Otto}, {\em The geometry of dissipative evolution equations: the porous
  medium equation}, Comm. Partial Differential Equations, 26 (2001),
  pp.~101--174.

\bibitem{PelZhang}
{\sc M.~Peletier and H.~Zhang}, {\em {Self-similar solutions of a fast
  diffusion equation that do not conserve mass.}}, Differ. Integral Equ., 8
  (1995), pp.~2045--2064.

\bibitem{MR1606339}
{\sc A.~Rodriguez, J.~L. Vazquez, and J.~R. Esteban}, {\em The maximal solution
  of the logarithmic fast diffusion equation in two space dimensions}, Adv.
  Differential Equations, 2 (1997), pp.~867--894.

\bibitem{SC1}
{\sc L.~Saloff-Coste}, {\em Precise estimates on the rate at which certain
  diffusions tend to equilibrium}, Math. Z., 217 (1994), pp.~641--677.

\bibitem{MR0463908}
{\sc G.~Talenti}, {\em Best constant in {S}obolev inequality}, Ann. Mat. Pura
  Appl. (4), 110 (1976), pp.~353--372.

\bibitem{MR1977429}
{\sc J.~L. V{\'a}zquez}, {\em Asymptotic behavior for the porous medium
  equation posed in the whole space}, J. Evol. Equ., 3 (2003), pp.~67--118.
\newblock Dedicated to Philippe B\'enilan.

\bibitem{VazBook}
\leavevmode\vrule height 2pt depth -1.6pt width 23pt, {\em The Porous Medium
  Equation.Mathematical Theory},  Oxford Mathematical Monographs. The Clarendon
Press, Oxford University Press, Oxford, 2007.

\bibitem{Vazquez2006}
\leavevmode\vrule height 2pt depth -1.6pt width 23pt, {\em Smoothing and decay
  estimates for nonlinear diffusion equations}, vol.~33 of Oxford Lecture Notes
  in Maths. and its Applications, Oxford Univ. Press, 2006.

\bibitem{MR1357953}
{\sc J.~L. V{\'a}zquez, J.~R. Esteban, and A.~Rodr{\'{\i}}guez}, {\em The fast
  diffusion equation with logarithmic nonlinearity and the evolution of
  conformal metrics in the plane}, Adv. Differential Equations, 1 (1996),
  pp.~21--50.

\end{thebibliography}
\end{document}